\documentclass[10pt]{article}
\usepackage[a4paper, margin=2cm]{geometry}

\usepackage[utf8]{inputenc}

\usepackage{amsmath, amssymb, amsfonts, mathtools, amsthm}
\usepackage{listings}
\usepackage{nccmath}
\usepackage{authblk}

\usepackage{bm} 

\usepackage{algorithm}
\usepackage{algpseudocode}
\usepackage{multirow}
\usepackage{booktabs}
\usepackage[]{mcode}

\newcommand*{\dittoclosing}{--- \raisebox{-0.5ex}{''} ---}

\newtheorem{theorem}{Theorem}[section]

\newtheorem{remark}[theorem]{Remark}

\usepackage{graphicx}
\usepackage{xcolor,url,hyperref}

\DeclareMathAlphabet{\mathbf}{OT1}{cmr}{bx}{it}
\newcommand{\vf}{\mathbf f}
\newcommand{\vg}{\mathbf g}
\newcommand{\vh}{\mathbf h}
\newcommand{\vr}{\mathbf r}
\newcommand{\vx}{\mathbf x}
\newcommand{\vu}{\mathbf u}
\newcommand{\vv}{\mathbf v}
\newcommand{\vw}{\mathbf w}

\newcommand{\alg}{\texttt{sketchAAA} }

\newcommand{\C}{{\mathbb C}}
\newcommand{\K}{{\mathbb K}}
\newcommand{\R}{{\mathbb R}}
\newcommand{\calA}{{\mathcal A}}
\newcommand{\vecop}{{\mathrm{vec}}}

\title{Randomized sketching of nonlinear~eigenvalue~problems}
\author[1]{Stefan Güttel}
\author[2]{Daniel Kressner}
\author[3]{Bart Vandereycken}
\affil[1]{\footnotesize Department of Mathematics, The University of Manchester, M139PL, Manchester, UK} 
\affil[2]{\footnotesize Institute of Mathematics, EPFL, CH-1015 Lausanne, Switzerland} 
\affil[3]{\footnotesize Section of Mathematics, University of Geneva, CH-1205 Geneva, Switzerland} 
\date{}

\begin{document}

\maketitle

\begin{abstract}
Rational approximation is a powerful tool to obtain accurate surrogates for nonlinear functions that are easy to evaluate and linearize. The interpolatory adaptive Antoulas--Anderson (AAA) method is one approach to construct such approximants numerically. For large-scale vector- and matrix-valued functions, however, the direct application of the set-valued variant of AAA becomes inefficient. We propose and analyze a new sketching approach for such functions called \alg that, with high probability, leads to much better approximants than previously suggested approaches while retaining efficiency. The sketching approach works in a black-box fashion where only evaluations of the nonlinear function at sampling points are needed. Numerical tests with nonlinear eigenvalue problems illustrate the efficacy of our approach, with speedups over 200 for sampling large-scale black-box functions without sacrificing on accuracy.
\end{abstract}

\section{Introduction}

This work is concerned with approximating vector-valued and matrix-valued functions, with a particular focus on  functions $F\colon \Sigma \to \mathbb C^{n\times n}$ that arise in the context of nonlinear eigenvalue problems (NLEVP) \begin{equation} \label{eq:nlevp}
F(z) \vv = 0, \quad \vv \not = 0.
\end{equation}

Recent algorithmic advances have made it possible to efficiently compute an accurate rational approximant of a scalar function $f\colon\Sigma \to \C$ on a compact (and usually discrete)  target set~$\Sigma$ in the complex plane~$\C$. 
Some of the available methods are the adaptive Antoulas--Anderson (AAA) algorithm~\cite{nakatsukasa2018aaa}, the rational Krylov fitting (RKFIT) algorithm~\cite{berljafa2017rkfit}, vector fitting~\cite{gustavsen1999rational}, minimal rational interpolation~\cite{Pradovera2020},
and methods based on L{\"o}wner matrices~\cite{mayo2007framework}. All of these methods can be used or adapted to approximate multiple scalar functions $f_1,\ldots,f_s$ on the target set $\Sigma$ simultaneously. In particular, there  are several variants of AAA that  approximate multiple functions by a family of rational interpolants sharing the same denominator, including the  set-valued AAA algorithm~\cite{lietaert2018automatic},  fastAAA~\cite{hochman2017fastaaa}, and weighted AAA~\cite{NGT20}. See also~\cite{EG21} for a discussion and comparison of some of these methods. 

 AAA-type algorithms can be used with very little user input and have enabled an almost black-box approximation of NLEVPs or transfer functions in model order reduction. However, the computation of $s$ degree-$d$ rational interpolants via the set-valued AAA algorithm involves computing the SVD of $d$ dense matrices of sizes~$s(|\Sigma| - k - 1)\times (k+1)$ for varying $k=1,2,\ldots,d$. Since these matrices differ only in $s$ rows and columns during the AAA iterations, the naive overall complexity of $O(s  |\Sigma| d^3)$ can be reduced to $O(s  |\Sigma| d^2 + d^4)$ with the SVD updating scheme in~\cite{lietaert2018automatic}. The greedy search of interpolation nodes in AAA also requires the repeated evaluation of the $s$~rational interpolants at all sampling points in~$\Sigma$ and the storage of the corresponding function values. As a consequence, the main use case for the multiple-function AAA variants to date have been problems that can be written in the split form 
\begin{equation}\label{eq:split form}
    F(z) = \sum_{i=1}^s f_i(z) A_i, 
\end{equation}
where $s$ is small (say, in the order of $s \approx 10$), $f_i$ are known scalar functions,  and  $A_i\in\mathbb{C}^{m\times n}$ are fixed coefficient matrices. While, in principle, it is always possible to write an arbitrary $m\times n$ matrix-valued function~$F$ in split form using $s = mn$ terms, it would be prohibitive to apply the set-valued AAA approach to large-scale problems in such a naive way. 

The work~\cite{EG19} suggested an alternative approach where the original (scalar-valued) AAA algorithm is applied to a scalar surrogate $f(z) = \vu^T F(z) \vv$ with random probing vectors $\vu$ and $\vv$, resulting in a rational interpolant in barycentric form
\begin{equation}\label{eq:rd}
    r^{(d)}(z) = \sum_{i=0}^d\frac{w_i \, f(z_i)}{z-z_i} \bigg/ \sum_{i=0}^d \frac{w_i}{z-z_i} 
\end{equation}
with support points~$z_i\in\C$ (the interpolation nodes) and weights~$w_i\in\C$.  Using this representation, a rational interpolant~$R^{(d)}$ of the original function $F$ is then  obtained by replacing  in \eqref{eq:rd} every occurrence of $f(z_i)$ by the evaluation $F(z_i)$:
\begin{equation}\label{eq:Rd}
     R^{(d)}(z) = \sum_{i=0}^d\frac{w_i \, F(z_i)}{z-z_i} \bigg/ \sum_{i=0}^d \frac{w_i}{z-z_i}.
\end{equation}
The intuition  is that both $F$ and $f$ will almost surely have the same region of analyticity, hence interpolating~$F$ using the same interpolation nodes and poles as for~$f$ should result in a good approximant. This surrogate approach indeed alleviates the complexity and memory issues even when $F$ has a large number of terms $s$ in its split form \eqref{eq:split form}, and it can also be applied if $F$ is only available as a black-box function returning evaluations $z\mapsto F(z)$. However, a comprehensive numerical comparison in~\cite{NGT20} in the context of solving NLEVPs~\eqref{eq:nlevp} has revealed that this surrogate approach is not always reliable and may  lead to poor accuracy. Indeed, the $\Sigma$-uniform error for the original problem
\begin{equation}
\| F-R^{(d)}\|_\Sigma := \max_{z\in\Sigma} \| F(z)-R^{(d)}(z)\|_F
\label{eq:sigerr}
\end{equation}
can be significantly larger than the error for the scalar surrogate problem $\| f - r^{(d)}\|_\Sigma$. {In order to mitigate this issue for black-box functions $F(z)$---i.e., those not available in split form \eqref{eq:split form}---a two-phase surrogate AAA algorithm with
cyclic Leja--Bagby refinement has been developed in~\cite{NGT20}. While this algorithm is indeed robust in the sense that it returns a rational approximant with user-specified accuracy, it is  computationally more expensive than the set-valued AAA algorithm and sometimes returns approximants of unnecessarily high degree; see, e.g., \cite[Table~5.2]{NGT20}.} 

In this work, we  propose and analyze a new sketching approach for matrix-valued functions  called \alg that, with high probability,  leads to much better approximants than the scalar surrogate AAA approach. 
At the same time, it remains equally efficient. While we demonstrate the benefits of the sketching  approach in combination with the set-valued AAA algorithm and mainly test it on functions from the updated NLEVP benchmark collection~\cite{betcke2013nlevp,higham2019updated}, the same technique can be combined with other linear interpolation techniques (including polynomial interpolation at adaptively chosen nodes). Our sketching approach is not limited to nonlinear eigenvalue problems either  and can be used for the approximation of any vector-valued function.  The key idea is to sketch  $\vf(z) = \vecop(F(z))$ using a thin (structured or unstructured) random probing matrix $V\in\C^{N\times \ell}$, i.e., computing samples of the form $V^T \vf(z_i)$, and then to apply the set-valued AAA algorithm to the $\ell$~components of the samples.  We provide a probabilistic assessment of the approximation quality of the resulting samples by building on  the  improved small-sample bounds of the matrix Frobenius norm in~\cite{gratton2018improved}.

The remainder of this work is organized as follows. In section~\ref{sec:alg} we briefly review the AAA algorithm~\cite{nakatsukasa2018aaa}  and the surrogate approach from \cite{EG19}, and then introduce our new probing forms using either tensorized or non-tensorized probing vectors. 
We also provide a pseudocode of the resulting \alg algorithm. 
Section~\ref{sec:analysis} is devoted to the  analysis of the approximants produced by  \alg\!\!.  Our main result for the case of non-tensorized probing, Theorem~\ref{thm:analysis}, provides an exponentially converging (in the number of samples $\ell$) probabilistic upper bound on the approximation error of the sketched problem compared to the approximation of the full problem. We also provide a weaker result for tensorized probing in the case of $\ell=1$, covering the original surrogate approach in \cite{EG19}. Several numerical tests in section~\ref{sec: experiments} on small to large-scale problems show that our new sketching approach is reliable and efficient. For some of the large-scale black-box problems we report speedup factors of over 200 compared to the  set-valued AAA algorithm implemented in \textsc{Matlab} while retaining comparable accuracy.

\section{Surrogate functions and the AAA algorithm}\label{sec:alg}

The AAA algorithm~\cite{nakatsukasa2018aaa} is a relatively simple but usually very effective method for obtaining a good rational approximant $r^{(d)}(z)\approx f(z)$ of the form \eqref{eq:rd}  for a given function $f \colon \C \supset \Omega  \to \C$. The approximation is sought on a finite discretization $\Sigma \subset \Omega$ and discretization sizes of $|\Sigma| \approx 10^3$--$10^4$ are not uncommon. The method iteratively selects support points (interpolation nodes) by adding to a previously computed set $\{z_0,z_1,\ldots, z_{d-1}\}$ of $d$ support points a new point $z_d\in\Sigma$ at which $\max_{z\in\Sigma} |r^{(d-1)}(z) - f(z)|$ is attained. After that, the rational approximant $r^{(d)}$ is calculated by computing $d+1$ weights $w_i \in \C$ as the minimizer of the linearized approximation error
\begin{equation}\label{eq:min L w}
 \min_{\vw\in\mathbb{C}^{d+1}} \| L\vw \|_2 \quad \text{such that} \quad \|\vw\|_2=1.
\end{equation}
Here, $L$ is a Löwner matrix of size $(|\Sigma|-d-1) \times  (d+1) $ defined as
\begin{equation}\label{eq:def L ji}
 L_{ij} = \frac{f(\hat z_i) - f(z_j)}{\hat z_i - z_j} 
\end{equation}
where $\{\hat z_1, \hat z_2, \ldots\} = \Sigma \setminus \{z_0, \ldots, z_d\}$. The minimization problem~\eqref{eq:min L w} can be solved exactly by SVD at a cost of $O(|\Sigma|d^2)$ flops. Since the matrix $L$ is only altered in one row and column after each iteration, updating strategies can be used to lower the computational cost of its SVD to $O(|\Sigma|d + d^3)$ flops; see, e.g.,~\cite{hochman2017fastaaa,lietaert2018automatic}. 
\textbf{}
It is not difficult to extend AAA for the approximation of a vector-valued function $\vf\colon \Omega \to \C^N$. Firstly, the vector-valued version $\vr^{(d)}$ of~\eqref{eq:rd} will map into $\C^N$ since the $\vf(z_i)$ are vector-valued. However, the support points $z_i$ and weights $w_i$ remain scalars. The selection of the support points can  still be done greedily by, at iteration~$d$, choosing a support point $z_d\in\Sigma$ that maximizes $\| \vr^{(d-1)}(z)  - \vf(z)\|$ on $\Sigma$ for some norm $\|\cdot\|$. In practice, the infinity or Euclidean norm usually work well but more care is sometimes needed when $\vf$ maps to different scales; see~\cite{lietaert2018automatic,NGT20}.  Likewise, the weights are computed from a block-Löwner matrix $L$ where each $L_{ij}$ in~\eqref{eq:def L ji} is now a column vector of length $N$ composed with $\vf(\hat z_i), \vf(z_j) \in \C^N$. The matrix $L$ is now of size $N(|\Sigma|-d-1) \times  (d+1) $, increasing the cost of its SVD to $O(N|\Sigma|d^2)$ flops. When computing these SVDs for degrees $1,2,\ldots,d$ as is required by AAA, the cumulative cost is $(N|\Sigma|d^3)$ when each SVD is computed independently. For large $N$, this  becomes prohibitive. Fortunately, the matrix $L$ is only changed in $s$ rows and columns during each iteration of AAA. One can therefore update the SVD as is done in~\cite{lietaert2018automatic}, reducing the overall complexity to $O(N|\Sigma|d^2 + d^4)$. However, this remains very costly for large $N$.

As explained in the introduction, a way to lower the computational cost for s vector-valued function $\vf$ is to work with a scalar surrogate function $g\colon \Omega \to \C$ that hopefully shares the same region of analyticity as $\vf$. In~\cite{EG19} this function was chosen with a tensorized probing vector:
\begin{equation}\label{eq:def g tensor}
g_\text{tens}(z) = (\vv \otimes \vu)^T \vf(z) \quad \text{with fixed random vectors $\vu,\vv \in \C^n$}.
\end{equation}
The reason for this construction with a tensor product is that~\cite{EG19} focused on nonlinear eigenvalue problems where $\vf(z) = \vecop(F(z))$ with $F(z) \in \C^{n \times n}$. This allows for an efficient evaluation $g_\text{tens}(z) = \vu^T F(z) \vv$, which becomes particularly advantageous when fast matrix-vector products with $F(z)$ are available. In the case that $F$ is in the split form~\eqref{eq:split form}, only a matrix-vector product with each $A_i$ is needed. A similar surrogate can be obtained for general $\vf\colon \Omega \to \C^N$ by using a full (non-tensorized) probing vector:
\begin{equation}\label{eq:def g notensor}
 g_\text{full}(z) = \vv^T \vf(z) \quad \text{with a fixed random vector $\vv\in\mathbb{C}^N$.}
\end{equation}
For both surrogate constructions, we apply AAA to $g_\text{tens}$ and $ g_\text{full}$ and use the computed support points $z_i$ and weights $w_i$ to define $\vr^{(d)}$ as in~\eqref{eq:rd},  replacing the scalar function values by the  vectors $\vf(z_i)$. Since the surrogate functions $g_\text{tens}$ are $ g_\text{full}$  are scalar-valued, the computational burden is clearly much lower than applying the set-valued AAA method to the vector-valued function~$\vf$.

\begin{figure}[h]
\begin{minipage}[b]{0.49\linewidth}
\centering
\includegraphics[width=\textwidth]{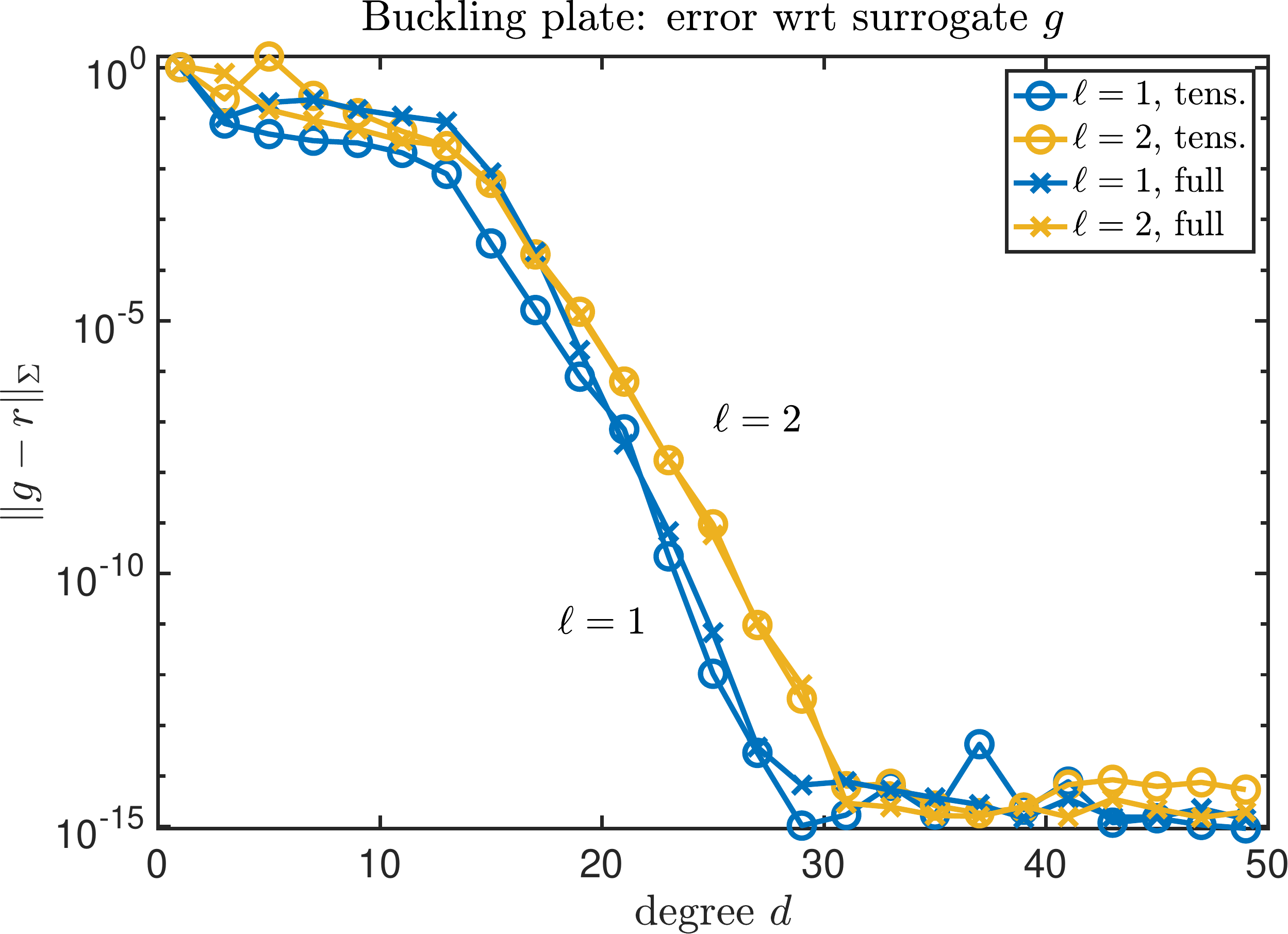}
\end{minipage}
\begin{minipage}[b]{0.49\linewidth}
\centering
\includegraphics[width=\textwidth]{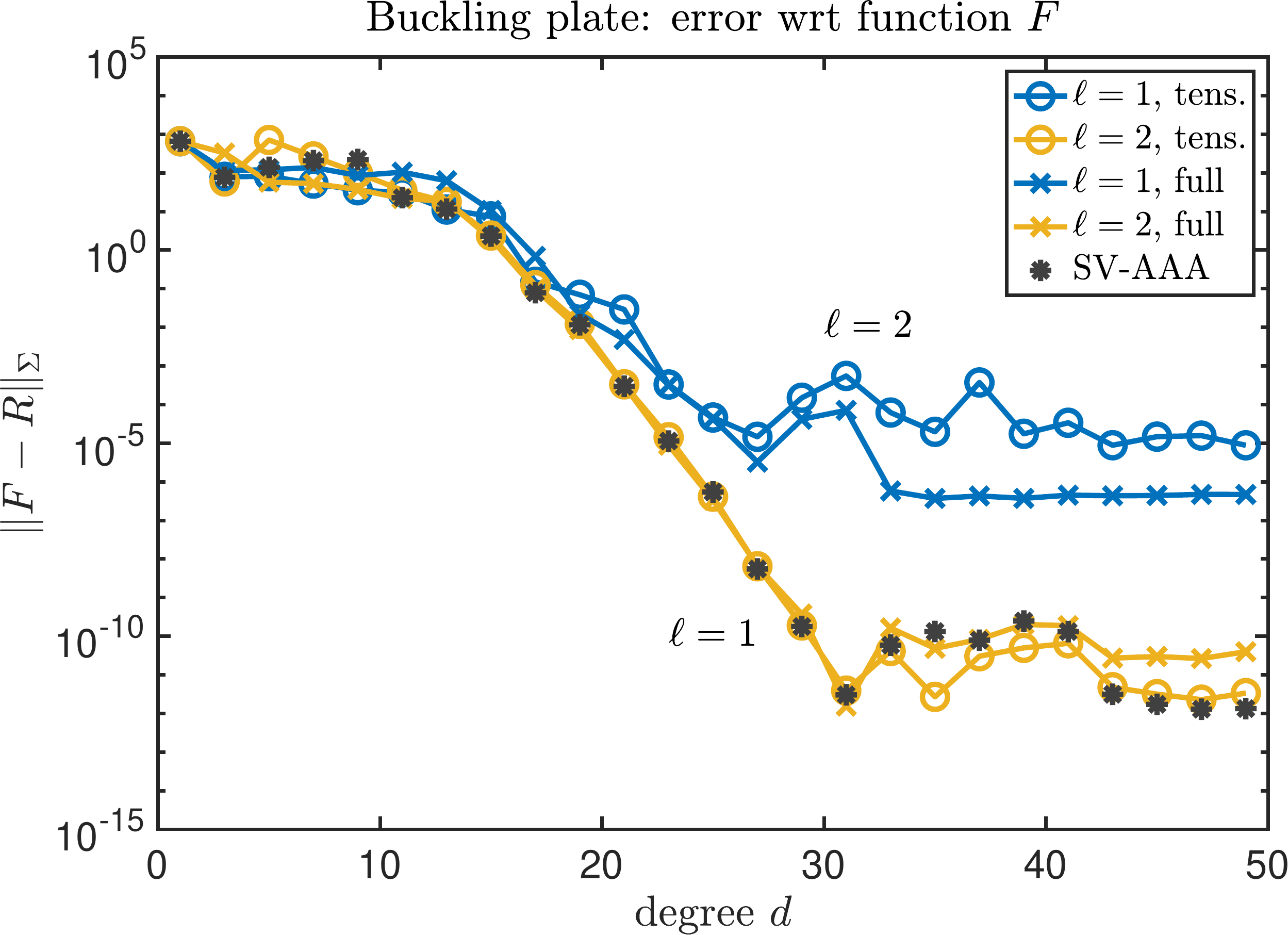}
\end{minipage}
\caption{Approximation errors of the rational interpolant for the \texttt{buckling\_plate} example as functions of the degree $d$. The error of the surrogate functions ($\ell=1$ or $\ell=2$ probing vectors, tensorized or not) is shown in the left panel, and on the right the error of the full interpolants is shown. While the errors of the surrogates all converge to machine precision (left), the error of the full interpolants stagnates when scalar surrogates ($\ell=1$) are used (right). With $\ell=2$ probing vectors the surrogate and full approximants converge similarly. On the right, the error of the set-valued AAA approximant applied the nine components of the original function~$F$ is indicated with dark grey dots.}
\label{fig:nice example}
\end{figure}

While computationally very attractive, the approach of building a scalar surrogate does unfortunately not always result in very accurate approximants. To illustrate, we consider the matrix-valued function~$F$ of the \texttt{buckling\_plate} example from~\cite{higham2019updated}. In section~\ref{sec:nlvep experiments} we will treat this problem among many others, and we refer to this section for the concrete experimental setup. In Figure~\ref{fig:nice example} the errors for both surrogates~\eqref{eq:def g tensor} and~\eqref{eq:def g notensor} are shown with the text label $\ell=1$. While AAA fully resolves each surrogate to an absolute accuracy well below $10^{-10}$, the absolute errors of the corresponding matrix-valued approximants to $F$ stagnate around~$10^{-5}$.

An important observation put forward in this paper is that by taking multiple random probing vectors (and thereby a vector-valued surrogate function $\vg$), the approximation error obtained with the set-valued AAA method can be improved, sometimes dramatically. In particular, we consider
\begin{equation}\label{eq:def g notensor vector valued}
 \vg_{\ell,\text{full}}(z) = V^T \vf(z) \quad \text{with a random matrix $V \in \C^{N \times \ell}$}
\end{equation}
and the tensorized version 
\begin{equation}\label{eq:def g tensor vector valued}
 \vg_{\ell,\text{tens}}(z) = [ \vv_1 \otimes \vu_1, \ldots,  \vv_\ell \otimes \vu_\ell]^T \vf(z) \quad \text{with random vectors $\vu_i,\vv_i \in \C^n$}.
\end{equation}
Both probing variants can be sped up if $F(z)$ (and hence $\vf(z)$) is available in the split form~\eqref{eq:split form} by precomputing the products of the random vectors with the matrices $A_1,\ldots,A_s$. The tensorized variant is computationally attractive when matrix vector products with $F(z)$ can be computed efficiently since $(\vv_i \otimes \vu_i)^T \vf(z)$ is the vectorization of $\vu_i^T F(z) \vv_i$. In both cases we obtain a vector-valued function with $\ell\ll N$ components to which the set-valued AAA method can readily be applied. The resulting algorithm \alg is summarized in Algorithm~\ref{alg:AAA-sketch}.

\begin{algorithm}[t]
\caption{AAA with random sketching (\texttt{sketchAAA})}\label{alg:AAA-sketch}
\begin{algorithmic}
\Require Nonlinear function $\vf\colon\mathbb{C} \supset \Sigma \to  \mathbb{C}^N$. Number of probing vectors $\ell$. 
Degree of approximation $d$.
\State 1.~Draw a random probing matrix $V \in \C^{N \times \ell}$.
\State 2.~Define (and possibly precompute on $\Sigma$) the surrogate $\vg_\ell(z) = V^T \vf(z)$.
\State 3.~Approximate $\vg_\ell$ on $\Sigma$ using the set-valued AAA method, resulting in a rational function in barycentric form with weights $\{w_i\}_{i=0}^d$ and support points $\{z_i\}_{i=0}^d$.
\State 4.~Return the rational approximant of $\vf$ in the form
\[
    \vr^{(d)}(z) = \sum_{i=0}^d\frac{w_i \, \vf(z_i)}{z-z_i} \bigg/ \sum_{i=0}^d \frac{w_i}{z-z_i} .
\]
\vspace*{-2mm}
\end{algorithmic}
\end{algorithm}

Sometimes using as few as $\ell=2$ probing vectors leads to very satisfactory results. This is indeed the case for our example in Figure~\ref{fig:nice example} where the approximation errors of both surrogates  \eqref{eq:def g notensor vector valued}  and \eqref{eq:def g tensor vector valued} are shown with the text label $\ell=2$. Both approximants converge rapidly to an absolute error below $10^{-10}$. Remarkably, the approximation error of the surrogates is essentially identical to the error of set-valued AAA approximant applied to the original function~$F$, which maps to $\C^{3 \times 3}$; see Figure~\ref{fig:nice example} (right).

\begin{figure}[ht]
\begin{minipage}[b]{0.49\linewidth}
\hspace*{-3mm}
\includegraphics[width=1.1\textwidth]{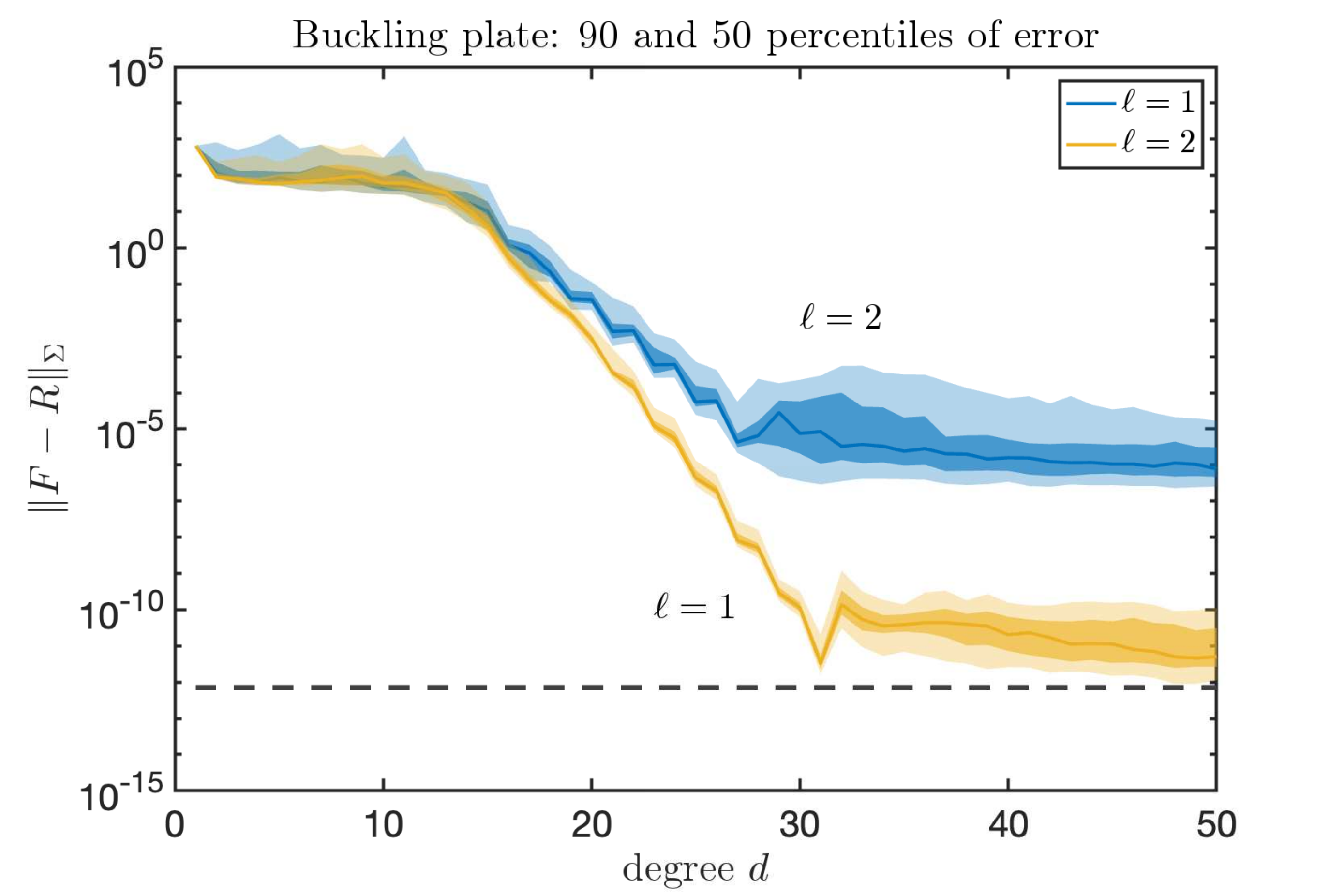}
\end{minipage}
\hfill
\begin{minipage}[b]{0.49\linewidth}
\hspace*{-3mm}
\includegraphics[width=1.1\textwidth]{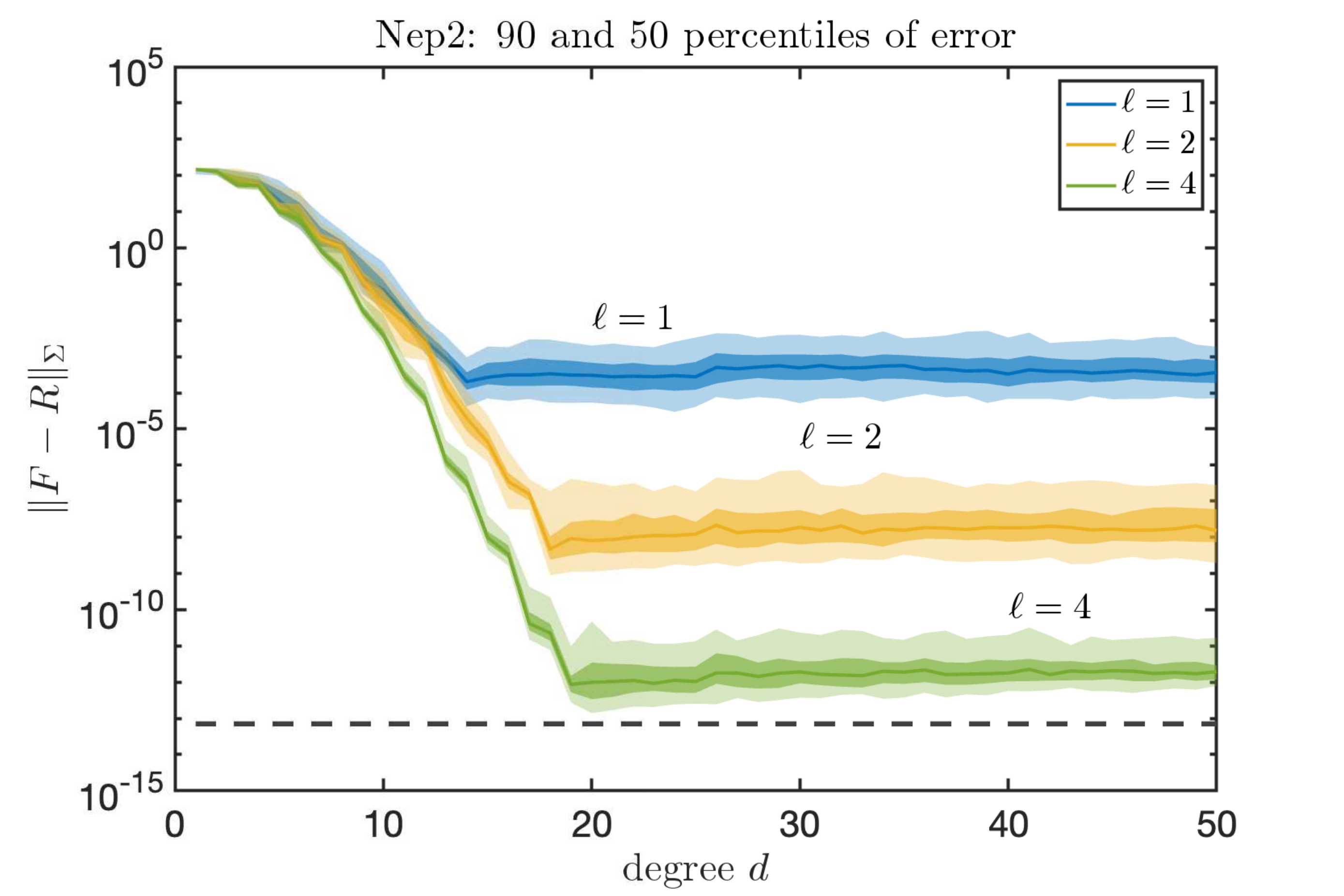}
\end{minipage}
\caption{Bands for the 90 and 50 percentiles of 100 random initializations of the surrogates (non-tensorized) for the \texttt{buckling\_plate} and \texttt{nep2} examples. The dashed line is the minimal error achieved by the  set-valued AAA method applied to the original function $F$.}
\label{fig:nice example pct}
\end{figure}

Since the surrogates are random, the resulting rational approximants are also random. Fortunately, their approximation errors concentrate well. This is clearly visible in Figure~\ref{fig:nice example pct} for the \texttt{bucking\_plate} and \texttt{nep2} examples from the NLEVP collection (see section~\ref{sec:nlvep experiments} for the concrete experimental setup), showing the  50-percentiles which are less than $10^{-1}$ wide. As a result, the random construction of the surrogates yields rational approximants that, with high probability, all have very similar approximation errors.  Figure~\ref{fig:nice example pct} also demonstrates that two probing vectors do not always suffice. For the \texttt{nep2} problem, $\ell=4$ probing vectors are required for a relative approximation error  close to machine precision.

Many more such examples are discussed in section~\ref{sec: experiments} and in almost all cases a relatively small number $\ell$ of probing vectors suffices to obtain accurate approximations. The next section provides some theoretical insights into this.

\section{Analysis of random sketching for function approximation}\label{sec:analysis}

In this section we analyze the effect of random sketching on the approximation error. More precisely, we show that for a given \emph{fixed} approximation of $\vf$ the corresponding approximation of the sketch $V^T \vf$ enjoys, with high probability, similar accuracy. Let us emphasize that this setting does not fully capture Algorithm~\ref{alg:AAA-sketch} because the approximation constructed by the algorithm depends on the random matrix $V$ in a highly nontrivial way. Nevertheless, we feel that the analysis explains the behavior of the algorithm for increasing $\ell$ and, more concretely, it allows us to cheaply and safely estimate the full error a posteriori via a surrogate (with an independent sketch).

\subsection{Preliminaries} \label{sec:prelims}

Let us consider a vector-valued function $\vh\colon \Omega \to \mathbb K^N$ on some domain $\Omega \subseteq \mathbb C$ and with $\mathbb K \in \{\mathbb R, \mathbb C\}$. Assuming $\vh \in L^2(\Omega, \K^N)$ we can equivalently view $\vh$ as an element of the tensor product space $\K^N \otimes L^2 (\Omega, \mathbb K)$ that induces the linear operator $\mathcal L_h\colon L^2 (\Omega, \mathbb K) \to \K^N$, $\mathcal L_h\colon v \mapsto \int_\Omega h(z) v(z)\, \mathrm{d}z$. Note that $\mathcal L_h$ is a Hilbert--Schmidt operator of rank at most~$N$. Applying the Schmidt decomposition~\cite[Theorem~4.137]{Hackbusch2019} implies the following result. 

\smallskip

\begin{theorem}  \label{thm:schmidt} 
Let $\vh \in L^2(\Omega, \K^N)$. Then there exist orthonormal vectors $\vu_1, \ldots, \vu_N \in \K^N$,  orthonormal functions $v_1,\ldots, v_N\in L^2(\Omega, \mathbb K)$, and scalars $\sigma_1 \ge \sigma_2 \ge \cdots \ge \sigma_N \ge 0$ such that
\[
 \vh(z) = \sum_{j = 1}^N \sigma_j \vu_j v_j(z).
\]
\end{theorem}

We note that the \emph{singular values} $\sigma_1,\ldots, \sigma_N$ are uniquely defined by $\vh$. 
By Theorem~\ref{thm:schmidt}, the norm of $\vh$ on $L^2(\Omega, \K^N)$ satisfies
\begin{equation}\label{eq:L2 norm}
\|\vh\|^2 := \int_\Omega \|\vh(z)\|_2^2\,\mathrm{d}z = \sigma_1^2 + \cdots + \sigma_N^2,
\end{equation}
where $\|\cdot\|_2$ denotes the Euclidean norm of a vector. Extending the corresponding notion for matrices, the \emph{stable rank} of $\vh$ is defined as $\rho := \|\vh\|^2 / \sigma_1^2$ and satisfies $1 \le \rho \le N$; see, e.g., \cite[\S2.1.15]{tropp2015introduction}.

Our analysis applies to an abstract approximation operator $\mathcal A: \mathcal W \to L^2(\Omega, \K^N)$ defined on some subspace $\mathcal W \subset L^2(\Omega, \K^N)$. 
We assume that $\mathcal A$ commutes with linear maps, that is,
\begin{equation} \label{eq:commuting}
\mathcal A( B \vf ) = B \mathcal A( \vf ), \quad \forall B \in \K^{N\times N}.
\end{equation}
This property holds for any approximation of the form
\[
\mathcal A( \vf ) = \sum_{i = 0}^d \vf(z_i) L_i(\,\cdot\,)    
\]
for \emph{fixed}  $z_j \in \Omega$ and $L_j \in L^2(\Omega, \K)$, with the tacit assumption that functions in $\mathcal W$ allow for point evaluations. In particular, this includes polynomial and rational interpolation  provided that the interpolation points and poles are considered fixed.
The relation~\eqref{eq:commuting} implies $\big(  {\widetilde \vg \atop 0} \big) = \mathcal A\big( {\widetilde \vf \atop 0} \big)$ for
$\widetilde \vf, \widetilde \vg: \Omega \to \mathbb C^m$, $m< N$. With a  slight abuse of notation we will simply write $\widetilde \vg = \mathcal A(\widetilde \vf)$.

We aim to analyze the difference between the full approximation error $\|\vf - \mathcal A(\vf) \|$ and the surrogate approximation error $\|  V^T \vf - \mathcal A( V^T \vf ) \|$. The assumption that $\mathcal A$ is considered fixed is satisfied by Algorithm~\ref{alg:AAA-sketch} \emph{if} another random probing matrix, independent from $V$, is used for the construction of the weights $\{w_i\}$ and the support points $\{z_i\}$. Our results then show when the surrogate approximation error provides a reliable and efficient a posteriori error estimate. Let us stress, however, that an actual implementation of Algorithm~\ref{alg:AAA-sketch} would normally use the same random matrix $V$ for constructing the approximation and measuring the error, creating nontrivial dependencies not captured by our analysis.

\subsection{Full (non-tensorized) probing vectors} \label{sec:fullanalysis}

The following theorem treats surrogates of the form~\eqref{eq:def g notensor vector valued}. It constitutes an extension of existing results~\cite{gratton2018improved,MR1337645} on small-sample matrix norm estimation.

\smallskip

\begin{theorem} \label{thm:analysis}
With the notation introduced in section~\ref{sec:prelims}, let
 $\rho$ denote the stable rank of $\vf - \mathcal A(\vf)$ for $\vf \in L^2(\Omega,\K^N)$ and let
$\widetilde V\in\K^{N\times \ell}$ be a real (for $\K = \R$) or complex (for $\K = \C$) Gaussian random matrix. Set $V =  \widetilde V / \sqrt{c \ell}$ with $c = 1$ for $\K = \R$ and $c = 2$ for $\K = \C$.
Then for any $\tau > 1$ we have 
\begin{equation} \label{eq:probunderestimate}
\mathrm{Prob}\Big( \|\vf - \mathcal A(\vf) \| \ge \tau \|  V^T \vf - \mathcal A( V^T \vf ) \|  \Big) \le \frac{1}{\Gamma(c\ell/2)} \gamma\big( c\ell/2, c \ell\rho \tau^{-2} / 2 \big),
\end{equation}
where $\gamma(s,x) := \int_0^x t^{s-1} e^{-t} \,\mathrm{d} t$ denotes the lower incomplete gamma function and
\begin{equation} \label{eq:proboverestimate}
\mathrm{Prob}\Big( \|\vf - \mathcal A(\vf) \| \le \tau^{-1} \|V^T \vf - \mathcal A( V^T \vf ) \|  \Big) \le \exp\Big(-\frac{c \ell}{2} \rho(\tau-1)^2\Big).
\end{equation}
\end{theorem}
\begin{proof} 
Applying the Schmidt decomposition from Theorem~\ref{thm:schmidt} to $\vh:= \vf - \mathcal A(\vf)$ and using~\eqref{eq:commuting} yields
\begin{equation} \label{eq:applysvd}
\|V^T \vf - \mathcal A( V^T \vf ) \|^2 = \|V^T \vh\|^2 = \sum_{j = 1}^N \sigma_j^2 \| V^T  \vu_j\|_2^2.
\end{equation}
Because $V$ is Gaussian and $\vu_1,\ldots, \vu_N$ is an orthonormal basis, it follows that $c \ell \| V^T \vu_j\|_2^2= \chi_j^2(c\ell)$ are mutually independent chi-squared variables with $c \ell$ degrees of freedom ($j = 1,\ldots,N$). Following well-established arguments~\cite{gratton2018improved,Roosta-K}, we obtain
\begin{eqnarray*}
\mathrm{Prob}( \|\vh\| \ge \tau \|V^T \vh\|) &=& 
\mathrm{Prob}\Big(c \ell \|V^T \vh\|^2 \le c \ell \tau^{-2} \|\vh\|^2 \Big) \\
&=&
\mathrm{Prob}\Big( \sum_{j = 1}^N \sigma_j^2 \chi_j^2(c \ell) \le c\ell \tau^{-2} \|\vh\|^2 \Big) \\
&\le & \mathrm{Prob}\big( \sigma_1^2 \chi_1^2(c \ell) \le c \ell \tau^{-2} \|\vh\|^2 \big) \\[2mm]
&=& \mathrm{Prob}\big( \chi_1^2(c \ell) \le c \ell \tau^{-2} \rho \big) \\
&=& \frac{1}{\Gamma(c\ell/2)} \gamma(c\ell/2,c \ell\rho \tau^{-2}/2),
\end{eqnarray*}
which proves~\eqref{eq:probunderestimate}.

The inequality~\eqref{eq:proboverestimate} follows directly from the proof of~\cite[Theorem 3.1]{gratton2018improved}, which establishes
\[
\mathrm{Prob}\Big( \frac{1}{c\ell} \sum_{j = 1}^N \sigma_j^2 \chi_j^2(c\ell) \ge \tau \|\vh\| \Big) \le \exp(-c\ell\rho(\tau-1)^2/2)
\]
and thus implies~\eqref{eq:proboverestimate}.\hfill
\end{proof}

\smallskip

To provide some intuition on the implications of Theorem~\ref{thm:analysis} for the complex case ($c = 2$), let us first note that
\[
  \Gamma(k) = (k-1)!,\quad   \gamma\big( k, \alpha \big) = \int_0^{\alpha} t^{k-1} e^{-t}\,\mathrm{d}t \approx \Big( \frac{\alpha}{2} \Big)^{k} \ \text{for}\  \alpha \approx 0.
 \]  
Setting $\alpha = 2\ell\rho \tau^{-2}$, this shows that the failure probability in~\eqref{eq:probunderestimate} is asymptotically proportional to $(e \rho \tau^{-2})^\ell$ and in turn, increasing $\ell$ will drastically reduce the failure probability provided that $\tau > \sqrt{e\rho}$. Specifically, for $\rho = 2$, $\tau = 10$, we obtain a failure probability of at most $2\%$ for $\ell = 1$ and $0.08\%$ for $\ell = 2$.
This means that if Algorithm~\ref{alg:AAA-sketch} returns an approximant that features a small error for a surrogate with $\ell = 2$ components, then the probability that the approximation error for the original function is more than ten times larger is below $0.08\%$. The probability that the error is more than hundred times larger is below $8\cdot 10^{-8}$.
On the other hand, if there exists a good approximant for $\vf$ then~\eqref{eq:proboverestimate} shows that it is almost guaranteed that the surrogate  function admits a nearly equally good approximant (which is hopefully found by the AAA algorithm). For the setting above, the probability that the error of the surrogate approximant is more than $3$ times larger than that of the approximant for the original function is less than $3\cdot 10^{-16}$.

\begin{remark} \label{remark:stablerank}
A large stable rank would lead to a less favourable bound~\eqref{eq:probunderestimate}, but there is good reason to believe that the stable rank of $\vf-\calA(\vf)$ remains modest in situations of interest. The algorithms discussed in this work are most meaningful when $\vf$ admits good rational approximants. More concretely, this occurs when $\epsilon_d:=\inf_{\vr_d} \|\vf-\vr_d\|$ decreases rapidly as  $d$ increases, where $\vr_d$ is a rational function of degree~$d$. In fact, $\epsilon_d$ decreases exponentially fast when $\vf$ is analytic in an open neighborhood of the target set $\Sigma$ (this is even true when the infimum is taken over the smaller set of  polynomials). As the $(d+1)th$ singular value of $\vf$ is bounded by $\epsilon_{d}$ this implies a rapid decay of the singular values and hence a small stable rank of $\vf$ and, likewise, of $\vf-\calA(\vf)$. The numerical experiment in section~\ref{sec: exp stats} confirms that the stable rank is indeed low.
\end{remark}

In practice, it may be convenient to use a real random matrix $V\in \R^{N\times \ell}$ for a complex-valued function $\vf \in L^2(\Omega,\K^N)$, for example, if the split form~\eqref{eq:split form} has real matrices $A_i$ and complex-valued $f_i$. Theorem~\ref{thm:analysis} extends to this situation by applying it separately to the real and imaginary parts of $\vf$ and using a union bound.

\subsection{Tensorized probing vectors}

The analysis of the surrogate~\eqref{eq:def g tensor vector valued} with tensorized probing vectors is significantly more complicated because tensorized random vectors fail to remain invariant under orthogonal transformations, an important ingredient in the proof of Theorem~\ref{thm:analysis}. As a consequence, the following partial results cover the case $\ell = 1$ only. They are direct extensions of analogous results for matrix norm estimation~\cite{BujanovicKressner2021}.

\smallskip

\begin{theorem} \label{theorem:rankone}
In the setting of Theorem~\ref{thm:analysis} with $\mathbb K = \mathbb R$, consider $V = \vv\otimes \vu$ for
real Gaussian random vectors $\vu \in \R^m$, $\vv \in \R^n$.
Then for any $\tau > 1$ we have 
\begin{equation} \label{eq:bound1rank1}
\mathrm{Prob}\big( \|\vf - \mathcal A(\vf) \| \ge \tau \sqrt{\rho} \|  V^T \vf - \mathcal A( V^T \vf ) \|  \big) \le \frac{2}{\pi} (2+\ln(1+2\tau)) \tau^{-1}
\end{equation}
and
\begin{equation}  \label{eq:bound2rank1}
\mathrm{Prob}\big( \|\vf - \mathcal A(\vf) \| \le \tau^{-1} \|V^T \vf - \mathcal A( V^T \vf ) \|  \big) \le \sqrt{2\tau} \exp(-\tau+2).
\end{equation}
\end{theorem}
\begin{proof}
Setting $\vh:= \vf - \mathcal A(\vf)$ and applying~\eqref{eq:applysvd} yields
\begin{eqnarray*}
\mathrm{Prob}( \|\vh\| \ge \tau \sqrt{\rho} \|  V^T \vh \| )
&=& \mathrm{Prob}\Big( \sigma_1^2 \ge \tau^2 \sum_{j = 1}^N \sigma_j^2 \|V^T  \vu_j \|_2^2 \Big) \\
&\le& \mathrm{Prob}\big( \|V^T  \vu_1 \|_2^2 \le \tau^{-2} \big).
\end{eqnarray*}
The last expression has been analyzed in the proof of Theorem~2.2 in~\cite{BujanovicKressner2021}, showing that it is bounded by the bound claimed in~\eqref{eq:bound1rank1}. 
Similarly, it follows from the proof of Theorem~2.4  in~\cite{BujanovicKressner2021} that the quantity
\[
\mathrm{Prob}\big( \|\vh \| \le \tau^{-1} \|V^T\vh\|  \big) = 
\mathrm{Prob}\Big( \tau^2 \sum_{j = 1}^N \sigma_j^2 \le \sum_{j = 1}^N \sigma_j^2 \|V^T  \vu_j \|_2^2 \Big)
\]
is bounded by the bound claimed in~\eqref{eq:bound2rank1}.\hfill
\end{proof}

\smallskip

While both failure probability bounds of Theorem~\ref{theorem:rankone} tend to zero as $\tau \to \infty$, the convergence predicted by~\eqref{eq:bound1rank1} is rather slow. It remains an open problem to establish better rates for $\ell > 1$.

\section{Applications and numerical experiments}\label{sec: experiments}

Algorithm~\ref{alg:AAA-sketch} was implemented in \textsc{Matlab}. The set-valued AAA (SV-AAA) we compare to (and also needed in Step~3 of Algorithm~\ref{alg:AAA-sketch}) is a modification of the original \textsc{Matlab} code\footnote{ \url{https://people.cs.kuleuven.be/karl.meerbergen/files/aaa/autoCORK.zip}}  from~\cite{lietaert2018automatic}. The SV-AAA code implements an updating strategy for the singular value decomposition of the L{\"o}wner matrix defined in~\eqref{eq:def L ji} to avoid its full recomputation when the degree is increased from $d$ to $d+1$. All default options in SV-AAA are preserved except that the expansion points for AAA are  greedily selected based on the maximum norm instead of the $2$-norm. 
In addition, the stopping condition is based on the relative maximum norm of the approximation that SV-AAA builds over the whole sampling set $\Sigma$. So, for example, if SV-AAA is applied to the scalar functions $f_1,\ldots,f_s$ with a  stopping tolerance $\texttt{reltol}$, then the algorithm terminates when the computed rational approximations $r_i^{(d)}\approx f_i$ satisfy
\[
 \frac{\max_{z\in\Sigma} \max_{i} | f_i(z)-r_i^{(d)}(z) |}{\max_{z\in\Sigma} \max_{i} | f_i(z) |} \leq \texttt{reltol}.
\]
Finally, we made minor code modifications that are up to 20\% faster \textsc{Matlab} code compared to the original version but still producing the same output. We also recall that Algorithm~\ref{alg:AAA-sketch} with $\ell=1$ and tensorised sketches is mathematically equivalent to the method proposed in~\cite{EG19}.

The experiments were run on an Intel i7-12\,700 with 64~GB~RAM. The software to reproduce our experiments is publicly available at~\cite{guttel_2024_11634529}.\footnote{See also \url{https://gitlab.unige.ch/Bart.Vandereycken/sketchAAA}}

\subsection{NLEVP benchmark collection} \label{sec:nlvep experiments}

We first test our approach for the non-polynomial problems in the NLEVP benchmark collection~\cite{betcke2013nlevp,higham2019updated} considered in~\cite{NGT20}. Table~\ref{table:nlevp} summarizes the key characteristics of these problems, including the matrix size~$n$ and the number~$s$ of terms in the split form~\eqref{eq:split form}. The target set is a disc or half disc specified in a meaningful way for each problem separately; see~\cite[Table~3]{NGT20} for details. 
We follow the procedure from~\cite{NGT20} for  generating the sampling points of the target set: 300 interior points are obtained by randomly perturbing a regular point grid covering a disc or half disc and another 100 points are placed uniformly on the boundary.  This gives a total of 400 points for the set  $\Sigma$. 

\begin{table}
    \centering
        \caption{Selected benchmark examples from the NLEVP collection.}
    \label{table:nlevp}
    \begin{tabular}{rlrr|rlrr} \toprule
         $\#$ & name & size & $s$ & $\#$ & Name & size & $s$  \\ \midrule
         1 & \verb!bent_beam! & 6 & 16         & 13 & \verb!pdde_symmetric! & 81 & 3 \\
         2 & \verb!buckling_plate! & 3 & 3     & 14 & \verb!photonic_crystal! & 288 &3 \\
         3 & \verb!canyon_particle! & 55 & 83  & 15 & \verb!pillbox_small! & 20 & 4 \\
         4 & \verb!clamped_beam_1d! & 100 &  3 & 16 & \verb!railtrack2_rep! & 1410 &3 \\
         5 & \verb!distributed_delay1! & 3 & 4 & 17 & \verb!railtrack_rep! & 1005 &3 \\
         6 & \verb!fiber! & 2400  &      3     & 18 & \verb!sandwich_beam! & 168 &3 \\
         7 & \verb!hadeler! & 200  &      3    & 19 & \verb!schrodinger_abc! & 10 & 4 \\
         8 & \verb!loaded_string! & 100 &  3   & 20 & \verb!square_root! & 20 & 2 \\
         9 & \verb!nep1! & 2  &       2        & 21 & \verb!time_delay! & 3 & 3 \\
         10 & \verb!nep2! & 3 &       10       & 22 & \verb!time_delay2! & 2 &3 \\
         11 & \verb!nep3! & 10 &      4        & 23 & \verb!time_delay3! & 10 & 3 \\
         12 & \verb!neuron_dde! & 2 & 5 & 
         24 &  \verb!gun! & 9956 & 4 \\ \bottomrule
    \end{tabular}
\end{table}

\subsubsection{Small NLEVPs}
We first focus on the small problems from the NLEVP collection, that is, the matrix size is below $1000$ (larger problems are considered separately in the following section).
Tables~\ref{tab:nlevp small low tol} and \ref{tab:nlevp small high tol} summarize the results with a stopping tolerance $\texttt{reltol}$ of $10^{-8}$ or $10^{-12}$, respectively. 

For each problem we show both the degree $d$ and the attained relative $\Sigma$-uniform approximation error 
\[
\texttt{relerr} =  \frac{\max_{z\in\Sigma} \max_{ij} | F_{ij}(z)-R_{ij}^{(d)}(z)|}{\max_{z\in\Sigma} \max_{ij} | F_{ij}(z)|}
\]
of four algorithmic variants:

\smallskip

\begin{description}
\item[\rm SV-AAA $f$] refers to the set-valued AAA algorithm applied to the $s$ scalar functions in the split form~\eqref{eq:split form} of each problem.
\item[\rm SV-AAA $F$] refers to the set-valued AAA algorithm applied to all entries of the matrix~$F$, which is only practical for relatively small problems.
\item[\rm \alg] is used with $\ell=1$ or $\ell=4$ and full (non-tensorized) probing vectors. The reported degrees and errors are averaged over 10 random realizations of probing vectors. 
\end{description}

\smallskip

We find that in all considered cases, the error achieved by  \alg with $\ell=4$ probing vectors is very close to the stopping tolerance. This is achieved without ever needing a degree significantly higher than that required by SV-AAA~$f$ and SV-AAA~$F$; an important practical consideration for solving NLEVPs (see section~\ref{sec:solvenlevp} below). No timings are reported here because all algorithms return their approximations within  a few milliseconds.

\begin{table}
\centering
 \caption{Small problems of the NLEVP collection. Stopping tolerance \texttt{reltol}  is $10^{-8}$. Sketching is done with full (non-tensorized) vectors.} \label{tab:nlevp small low tol}
\begin{tabular}{r|cccccccc} 
$\#$  & \multicolumn{2}{c}{ SV-AAA $f$ } & \multicolumn{2}{c}{ SV-AAA $F$ }  & \multicolumn{2}{c}{ \alg $\ell=1$ } & \multicolumn{2}{c}{ \alg $\ell=4$ } \\
& degree & \texttt{relerr}  & degree & \texttt{relerr} & degree & \texttt{relerr} &  degree & \texttt{relerr}  \\
\hline
     1  &       9  &  4.5e-10  &   8  &  2.9e-09  &  5.7  &  1.3e-04   &   8  &  1.8e-08 \\  
      2  &      25  &  3.4e-09  &  25  &  3.4e-09  &  22.9  &  1.3e-05   &  25  &  3.5e-09 \\  
      3  &      16  &  7.9e-11  &  14  &  1.6e-09  &  9.7  &  1.3e-05   &  12.3  &  7.2e-08 \\  
      4  &      12  &  8.7e-09  &  12  &  6.6e-09  &  12  &  6.0e-06   &  12.2  &  6.4e-09 \\  
      5  &       8  &  1.8e-10  &   8  &  9.5e-11  &   7  &  5.5e-05   &  7.8  &  2.3e-09 \\  
      7  &       4  &  7.4e-09  &   4  &  7.1e-09  &  3.8  &  8.6e-07   &  4.2  &  2.3e-08 \\  
      8  &       3  &  3.2e-16  &   3  &  3.0e-16  &   3  &  5.8e-14   &   3  &  2.4e-15 \\  
      9  &      22  &  9.4e-09  &  21  &  7.7e-09  &  21  &  7.7e-09   &  21  &  7.7e-09 \\  
      10  &      15  &  1.7e-09  &  13  &  7.1e-09  &  10.7  &  2.0e-04   &  13  &  2.4e-08 \\  
      11  &      10  &  1.4e-09  &   9  &  3.8e-09  &   9  &  3.7e-06   &   9  &  8.4e-09 \\  
      12  &      15  &  9.8e-09  &  14  &  5.5e-09  &  14  &  2.6e-07   &  14  &  6.2e-09 \\  
      13  &       9  &  2.9e-09  &   9  &  3.1e-09  &  8.2  &  2.5e-05   &   9  &  4.0e-09 \\  
      14  &       7  &  5.7e-11  &   6  &  3.8e-11  &  6.3  &  1.4e-08   &  5.3  &  9.1e-09 \\  
      15  &       9  &  7.0e-11  &   7  &  3.8e-09  &  5.9  &  3.4e-05   &   7  &  3.8e-09 \\  
      18  &      12  &  1.1e-14  &   6  &  1.6e-09  &  5.8  &  6.9e-08   &  5.7  &  4.1e-09 \\  
      19  &      13  &  1.2e-09  &  12  &  4.3e-09  &  11  &  1.0e-05   &  12  &  6.5e-09 \\  
      20  &      11  &  2.2e-09  &  11  &  2.8e-09  &  10.8  &  2.9e-08   &  10.9  &  6.2e-09 \\  
      21  &      14  &  6.3e-09  &  14  &  1.3e-09  &  14  &  1.5e-09   &  14  &  1.3e-09 \\  
      22  &      14  &  5.9e-09  &  14  &  4.7e-09  &  14  &  1.1e-08   &  14  &  4.6e-09 \\  
      23  &      18  &  6.1e-09  &  14  &  5.0e-09  &  14  &  1.0e-08   &  14  &  4.4e-09 \\  
 \end{tabular}
 \end{table}

\begin{table}
\centering
 \caption{Small problems of the NLEVP collection. Stopping tolerance \texttt{reltol}  is $10^{-12}$. Sketching with full (non-tensorized) vectors.}\label{tab:nlevp small high tol}
\begin{tabular}{l|cccccccc} 
$\#$  & \multicolumn{2}{c}{ SV AAA $f$ } & \multicolumn{2}{c}{ SV AAA $F$ }  & \multicolumn{2}{c}{ \alg $\ell=1$ } & \multicolumn{2}{c}{ \alg $\ell=4$ } \\
& degree & \texttt{relerr}  & degree & \texttt{relerr} & degree & \texttt{relerr} &  degree & \texttt{relerr}  \\\hline
      1  &      12  &  4.9e-14  &  11  &  4.8e-13  &  7.8  &  1.7e-06   &  10.8  &  1.8e-11 \\  
      2  &      30  &  5.6e-13  &  30  &  5.6e-13  &  26.8  &  2.2e-07   &  30  &  8.3e-13 \\  
      3  &      23  &  1.2e-14  &  21  &  2.4e-13  &  15.3  &  1.1e-07   &  19.4  &  7.2e-12 \\  
      4  &      15  &  4.9e-13  &  15  &  5.9e-13  &  15  &  3.6e-08   &  15  &  5.8e-13 \\  
      5  &      10  &  3.8e-15  &   9  &  4.0e-13  &   9  &  2.9e-07   &   9  &  5.4e-13 \\  
      7  &      10  &  3.0e-13  &  10  &  3.3e-13  &   9  &  1.0e-08   &  10  &  4.9e-13 \\  
      8  &       3  &  3.2e-16  &   3  &  3.0e-16  &   3  &  5.8e-14   &   3  &  2.4e-15 \\  
      9  &      27  &  3.5e-13  &  27  &  2.9e-13  &  27  &  3.0e-13   &  27  &  3.0e-13 \\  
      10  &      18  &  9.5e-14  &  17  &  1.5e-13  &  12.6  &  2.1e-05   &  17  &  1.3e-12 \\  
      11  &      13  &  9.2e-14  &  12  &  8.6e-14  &  11  &  3.5e-08   &  12  &  1.2e-13 \\  
      12  &      19  &  8.0e-14  &  18  &  6.6e-14  &  17  &  6.5e-09   &  18  &  7.7e-14 \\  
      13  &      11  &  7.2e-13  &  11  &  6.6e-13  &  10.6  &  1.9e-07   &  11  &  5.7e-13 \\  
      14  &       8  &  9.4e-16  &   7  &  6.7e-16  &  7.8  &  8.8e-12   &  7.3  &  2.4e-15 \\  
      15  &      13  &  4.6e-15  &  11  &  2.0e-13  &  9.3  &  6.6e-08   &  11  &  1.7e-13 \\  
      18  &      16  &  7.6e-16  &  11  &  2.1e-13  &  10.8  &  6.0e-10   &  10.6  &  5.1e-13 \\  
      19  &      15  &  2.8e-13  &  15  &  5.3e-14  &  13  &  3.2e-07   &  14.8  &  8.9e-13 \\  
      20  &      15  &  3.7e-13  &  15  &  5.5e-13  &  14.9  &  3.1e-12   &  15.1  &  8.8e-13 \\  
      21  &      18  &  8.0e-14  &  18  &  4.6e-14  &  17  &  7.9e-12   &  17.8  &  2.4e-13 \\  
      22  &      18  &  8.0e-14  &  18  &  8.6e-14  &  17  &  4.3e-10   &  18  &  7.7e-14 \\  
      23  &      23  &  1.6e-13  &  18  &  2.3e-13  &  17  &  2.5e-10   &  18  &  2.6e-13 \\    
 \end{tabular}
 \end{table}

\subsubsection{Large NLEVPs} 
 
This section considers the large problems from the NLEVP collection listed in Table~\ref{table:nlevp}, all with problem sizes above 1000. More precisely, we consider problem 6 (size 2400), 16 (size 1410), 17 (size 1005), and 24 (size~9956). For such problem sizes, the application of set-valued AAA to each component of $F$ is no longer feasible and hence this algorithm is not reported. In order to simulate a truly black-box sampling of these eigenvalue problems when using full (non-tensorized) probing vectors as in~\eqref{eq:def g notensor vector valued}, we use the MATLAB function shown in Algorithm~\ref{alg:sparseprobe}. This function obtains the samples $V^T \vecop(F(z_i))$ without forming the large $n^2\times \ell$ Gaussian random matrix $V$ explicitly. Instead, the sparsity pattern of $F(z)$ is inferred on-the-fly as the sampling proceeds. In the case of tensorized probing as in \eqref{eq:def g tensor vector valued}, we can exploit sparsity more easily by computing $W = F(z_i)[\vv_1,\ldots,\vv_\ell]$, followed by the computation of  $\vu_1^T \vw_1, \ldots,\vu_\ell^T \vw_\ell$ where $\vw_i$ is the $i$th column of $W$.

The execution times are now more significant and reported in Tables~\ref{tab:nlevp large dense} and~\ref{tab:nlevp large tensor}, together with the required degree~$d$ and the achieved relative approximation error.  Table~\ref{tab:nlevp large tensor} shows that tensorized sketches lead to a faster algorithm with similar accuracy compared to the full sketches reported in Table~\ref{tab:nlevp large dense}.
Like for the small NLEVP problems in the previous section, we find that \alg with $\ell=4$ probing vectors is reliable and yields an approximation error close to the stopping tolerance~\texttt{reltol}.  

We note that these four problems are also available in split form \eqref{eq:split form} and both the non-tensorized and tensorized probing can be further sped up (sometimes significantly) by precomputing the products of the random vectors with the matrices $A_1,\ldots,A_s$. This is particularly the case  for the \texttt{gun} problem number~24 for which \alg spends most of its time on the evaluation of $F(z_i)$ at the support points~$z_i$. Exploiting the split form reduces this time drastically, which can be seen from the rows in Tables~\ref{tab:nlevp large dense} and \ref{tab:nlevp large tensor} labelled with ``24*''. The case $\ell=4$ is particularly interesting as, coincidentally, the problem also has $s=4$ terms and, in turn, the set-valued approximations for SV-AAA~$f$ and \alg both involve four functions. For tensorized sketching (Table~\ref{tab:nlevp large tensor}), \alg is faster than SV-AAA~$f$ while returning an accurate approximation of lower degree. This nicely demonstrates that our approach to exploiting the split form is genuinely different from the set-valued AAA algorithm in~\cite{lietaert2018automatic}: \alg takes the contributions of the coefficient matrices $A_i$ into account while the set-valued AAA algorithm only approximates the scalar functions~$f_i$ and is blind to $A_i$. Section~\ref{sec:split} below illustrates this further. However, it should be noted that for problems 6 and 24, the approximations computed by \alg ($\ell=4$) achieve an error very close to the targeted \texttt{reltol}, while SV-AAA~$f$ over-resolves. (Problems 16 and 17 are rational eigenvalue problems and hence they are resolved close to machine precision by all methods.)

\begin{algorithm}[t]
\vspace*{-2mm}
\begin{lstlisting}
function [vals, ind] = sparseprobe(F, Z, ell)
vals = zeros(ell,length(Z)); ind = []; V = []; 
for i = 1:length(Z)
    Fz = F(Z(i)); indi = find(Fz);             % evaluate F & find nonzeros
    if ~isequal(indi,ind)                      % reduce use of slow setdiff
        sd = setdiff(indi,ind);          
        if ~isempty(sd)                        % any new nonzeros?
            ind = [ind; sd];                   % add them to ind
            V = [V, randn(ell,length(sd))];    % and expand V
            [ind, r] = sort(ind); V = V(:,r);  % preserve index order
        end
    end
    vals(:,i) = V*Fz(ind);                     % probe nonzeros
end
\end{lstlisting}
\caption{MATLAB code for probing a sparse $F(z)$ on the target set \texttt{Z} using \texttt{ell}  random vectors. Returns samples in \texttt{vals} and the indices of $F$'s nonzeros in~\texttt{ind}. \label{alg:sparseprobe}}
\end{algorithm}

  \begin{table}
\centering
 \caption{Large problems of the NLEVP collection.  Sketching with full (non-tensorized) vectors. The~split form is not exploited by \alg except for problem 24*.} \label{tab:nlevp large dense}
\begin{tabular}{llcccccc} \toprule
\texttt{reltol} & $\#$  & \multicolumn{2}{c}{ SV-AAA $f$ } &  \multicolumn{2}{c}{ \alg $\ell=1$ } & \multicolumn{2}{c}{ \alg $\ell=4$ } \\
  & & degree & error  & degree & error &  degree & error  \\
 & & \multicolumn{2}{c}{ time (s.) }  & \multicolumn{2}{c}{ time (s.) } &  \multicolumn{2}{c}{ time (s.)} \\\midrule
 1e-08   
     & 6  &      14  &  9.0e-11  &  8.8  &  2.7e-05   &  9.8  &  2.7e-07 \\  
     &   &     \multicolumn{2}{c}{ 0.014 }   &  \multicolumn{2}{c}{ 0.039 }   &  \multicolumn{2}{c}{ 0.045 } \\  
     & 16  &       3  &  5.5e-16  &   3  &  6.4e-14   &  3.7  &  1.4e-14 \\  
     &   &     \multicolumn{2}{c}{ 0.006 }   &  \multicolumn{2}{c}{ 0.280 }   &  \multicolumn{2}{c}{ 0.299 } \\  
     & 17  &       3  &  3.8e-16  &  3.9  &  3.0e-13   &  3.7  &  2.1e-14 \\  
     &   &     \multicolumn{2}{c}{ 0.003 }   &  \multicolumn{2}{c}{ 0.461 }   &  \multicolumn{2}{c}{ 0.482 } \\  
     & 24  &      11  &  2.9e-11  &  7.3  &  1.5e-05   &  8.2  &  1.9e-08 \\  
     &   &     \multicolumn{2}{c}{ 0.016 }   &  \multicolumn{2}{c}{ 1.018 }   &  \multicolumn{2}{c}{ 1.097 } \\[1.5mm]
     & 24*  &    \multicolumn{2}{c}{\dittoclosing}   &  \multicolumn{2}{c}{ 0.019 }   &  \multicolumn{2}{c}{ 0.024 } \\ 
\midrule
1e-12   
     & 6  &      19  &  1.5e-14  &  14.3  &  3.6e-07   &  14.9  &  2.5e-11 \\  
     &   &     \multicolumn{2}{c}{ 0.006 }   &  \multicolumn{2}{c}{ 0.041 }   &  \multicolumn{2}{c}{ 0.045 } \\  
     & 16  &       3  &  5.5e-16  &   3  &  6.4e-14   &  3.7  &  1.4e-14 \\  
     &   &     \multicolumn{2}{c}{ 0.007 }   &  \multicolumn{2}{c}{ 0.281 }   &  \multicolumn{2}{c}{ 0.297 } \\  
     & 17  &       3  &  3.8e-16  &  3.9  &  3.0e-13   &  3.7  &  2.1e-14 \\  
     
     &   &     \multicolumn{2}{c}{ 0.002 }   &  \multicolumn{2}{c}{ 0.458 }   &  \multicolumn{2}{c}{ 0.481 } \\  
     & 24  &      15  &  3.5e-15  &  11  &  8.0e-08   &  12.8  &  2.6e-12 \\  
     &   &     \multicolumn{2}{c}{ 0.020 }   &  \multicolumn{2}{c}{ 1.024 }   &  \multicolumn{2}{c}{ 1.100 } \\[1.5mm]
     & 24*  &   \multicolumn{2}{c}{\dittoclosing}      &  \multicolumn{2}{c}{ 0.022 }   &  \multicolumn{2}{c}{ 0.030 } \\       \bottomrule  
 \end{tabular}
 \end{table}
 
   \begin{table}
\centering
 \caption{Large problems of the NLEVP collection.  Sketching with tensorized vectors. The split form is not exploited by \alg except for problem 24*. For $\ell=1$, the method is equivalent to that from~\cite{EG19}.} \label{tab:nlevp large tensor}
\begin{tabular}{llcccccc} \toprule
\texttt{reltol} & $\#$  & \multicolumn{2}{c}{ SV-AAA $f$ } &  \multicolumn{2}{c}{ \alg $\ell=1$ } & \multicolumn{2}{c}{ \alg $\ell=4$ } \\
  & & degree & error  & degree & error &  degree & error  \\
 & & \multicolumn{2}{c}{ time (s.) }  & \multicolumn{2}{c}{ time (s.) } &  \multicolumn{2}{c}{ time (s.)} \\\midrule
1e-08   
     & 6  &      14  &  9.0e-11  &  9.4  &  2.0e-05   &  9.8  &  2.9e-07 \\  
     &   &     \multicolumn{2}{c}{ 0.014 }   &  \multicolumn{2}{c}{ 0.016 }   &  \multicolumn{2}{c}{ 0.027 } \\  
     & 16  &       3  &  5.5e-16  &  3.1  &  2.6e-14   &  3.9  &  1.5e-15 \\  
     &   &     \multicolumn{2}{c}{ 0.006 }   &  \multicolumn{2}{c}{ 0.076 }   &  \multicolumn{2}{c}{ 0.106 } \\  
     & 17  &       3  &  3.8e-16  &  3.6  &  1.2e-14   &  3.5  &  3.1e-15 \\  
     &   &     \multicolumn{2}{c}{ 0.003 }   &  \multicolumn{2}{c}{ 0.098 }   &  \multicolumn{2}{c}{ 0.135 } \\  
     & 24  &      11  &  2.9e-11  &  6.9  &  1.4e-05   &  8.1  &  2.0e-08 \\  
     &   &     \multicolumn{2}{c}{ 0.016 }   &  \multicolumn{2}{c}{ 0.473 }   &  \multicolumn{2}{c}{ 0.679 } \\[1.5mm]
     & 24*  &     \multicolumn{2}{c}{\dittoclosing}   &  \multicolumn{2}{c}{ 0.009 }   &  \multicolumn{2}{c}{ 0.011 } \\  
\midrule
1e-12   
     & 6  &      19  &  1.5e-14  &  14.1  &  3.9e-07   &  15  &  2.6e-11 \\  
     &   &     \multicolumn{2}{c}{ 0.007 }   &  \multicolumn{2}{c}{ 0.017 }   &  \multicolumn{2}{c}{ 0.028 } \\  
     & 16  &       3  &  5.5e-16  &  3.1  &  2.6e-14   &  3.9  &  1.5e-15 \\  
     &   &     \multicolumn{2}{c}{ 0.007 }   &  \multicolumn{2}{c}{ 0.076 }   &  \multicolumn{2}{c}{ 0.108 } \\  
     & 17  &       3  &  3.8e-16  &  3.6  &  1.2e-14   &  3.5  &  3.1e-15 \\  
     &   &     \multicolumn{2}{c}{ 0.002 }   &  \multicolumn{2}{c}{ 0.097 }   &  \multicolumn{2}{c}{ 0.133 } \\  
     & 24  &      15  &  3.5e-15  &  10.4  &  3.5e-07   &  13  &  4.1e-13 \\  
     &   &     \multicolumn{2}{c}{ 0.020 }   &  \multicolumn{2}{c}{ 0.477 }   &  \multicolumn{2}{c}{ 0.673 } \\[1.5mm]
     & 24*  &   \multicolumn{2}{c}{\dittoclosing}     &  \multicolumn{2}{c}{ 0.013 }   &  \multicolumn{2}{c}{ 0.018 } \\  
\bottomrule  
 \end{tabular}
 \end{table}

\subsection{An artificial example on the split form}\label{sec:split}  The difference in the approximation degrees returned by SV-AAA~$f$ and \alg becomes particularly pronounced when cancellations occur between the different terms of the split form or when the coefficients are of significantly different scales. To demonstrate this effect by an extreme example, consider 
$$
F(z) = |z| \cdot 10^{-8} B + \sin(\pi z) C,
$$
where $z \in [-1,1]$ and $B,C \in \R^{10 \times 10}$ are random matrices of unit spectral norm. For the split form, we simply take the functions $f_1(z) = |z|$ and $f_2(z) = \sin(\pi z)$. The sampling set $\Sigma$ contains 100 equidistant points in $[-1,1]$. Since SV-AAA $f$ scales the functions $f_i$ to have unit $\infty$-norm on $\Sigma$, the results would remain the same if we took $f_1(z) = 10^{-8} |z|$ and $f_2(z) = \sin(\pi z)$. In Table~\ref{tab:artificial cancellation}, we clearly see that SV-AAA $f$ overestimates the degree since it puts too much emphasis on resolving $f_1(z)$.

  \begin{table}
\centering
 \caption{Artificial example demonstrating scaling issues arising with the split form used by SV-AAA~$f$. Sketching with full (non-tensorized) vectors. Our method \alg is immune to such effects.} \label{tab:artificial cancellation}
\begin{center}
\begin{tabular}{lcccccc} \toprule
\texttt{reltol}   & \multicolumn{2}{c}{ SV-AAA $f$ } &  \multicolumn{2}{c}{ SV-AAA $F$ } & \multicolumn{2}{c}{ \alg $\ell=4$ } \\
  &  degree & error  & degree & error &  degree & error   \\\midrule
 1e-08   
        &      24  &  1.1e-09     &  8  &  4.1e-09   &  8 &  4.1e-09 \\
1e-12   
       &      29  &  5.0e-13  &  18  &  2.5e-13   &  18 &  2.9e-13 \\   \bottomrule
 \end{tabular}
 \end{center}
 \end{table}

\subsubsection{Impact on numerical solution of NLEVP} \label{sec:solvenlevp}

As explained in \cite[Section~6]{GT17} and further analyzed in \cite[Section~2]{NGT20}, an accurate uniform rational approximation $R^{(d)}\approx F$ on the target set is crucial for a robust  linearization-based NLEVP solver. We refer in particular to \cite[eq.~(2.2)]{NGT20} and its discussion, which argues that if $\| F- R^{(d)}\|_\Sigma \leq \varepsilon \| F \|_\Sigma$ on a sufficiently fine discretization $\Sigma$ of a compact set $\Omega \subset \Sigma$, and $F,R^{(d)}$ are continuous on $\Omega$, then any eigenpair $(\lambda,v)$ of $R^{(d)}$, $\lambda\in \Omega$, will have a small backward error as an eigenpair of~$F$. Conversely, if $\mu\in \Omega$ is not an eigenvalue of $F$ (i.e., $F(\mu)$  is nonsingular), then a sufficiently accurate approximant $R^{(d)}$ is also nonsingular at~$\mu$; see \cite[Section~6]{GT17}.

Once an accurate rational approximant $R^{(d)}\approx F$ is obtained, it can be linearized in various ways; see, e.g., \cite{lietaert2018automatic,EG19,NGT20}. Specifically,  \cite[Theorem~3]{EG19} derives a (strong) linearization $L(z) = A - zB$ with $A,B\in \mathbb K^{dn\times dn}$ of a rational matrix function $R^{(d)}(z)$ in barycentric form~\eqref{eq:Rd}, as returned by \alg\!. We include this theorem here for completeness and notational consistency.

\smallskip

 \begin{theorem}\label{thm:nleigs}
 \label{thm:linearization} 
 Given an $n\times n$ rational function $R^{(d)}(z)$ in barycentric form~\eqref{eq:Rd} with weights $w_j$ and support points $z_j$. Let  $h_j,k_j,h_j$ be finite parameters so that 
 \[
     \beta_j k_j = -w_{j-1}/w_j, \quad \beta_j h_j = -z_jw_{j-1}/w_j \quad \text{for} \ \ j=1,\ldots,d.
 \] 
 Then the $nd\times nd$  pencil $L(z)=A-zB$ with
 \begin{align*}
 A &= \begin{small}\begin{bmatrix}
 h_d F(z_0) & h_d F(z_1) & \cdots & h_{d} F(z_{d-2}) & h_{d}F(z_{d-1})-h_d z_{d-1} F(z_{d})/\beta_d \\
  h_1 z_0 I & h_1\beta_1 I \\
  & \ddots & \ddots \\
  & & h_{d-2} z_{d-3} I & h_{d-2}\beta_{d-2} I & \\
  & & & h_{d-1} z_{d-2} I & h_{d-1}\beta_{d-1}I
 \end{bmatrix},\end{small} 
 \\
 B &= \begin{small}\begin{bmatrix} 
 k_d F(z_0) & k_d F(z_1) & \cdots & k_d F(z_{d-2}) & k_d F(z_{d-1}) - h_d F(z_d)/\beta_d \\
  h_1I & k_1\beta_1 I \\
  & \ddots & \ddots \\
  & & h_{d-2}I & k_{d-2}\beta_{d-2} I & \\
  & & & h_{d-1}I & k_{d-1}\beta_{d-1} I
 \end{bmatrix},\end{small} 
 \end{align*}
 is a strong linearization of $R^{(d)}(z)$.
 \end{theorem}

\smallskip

It follows from this theorem that the eigenvalues of $R^{(d)}(z)$ can be computed by solving a generalized eigenvalue problem with $(A,B)$, e.g., iteratively by applying a rational Krylov subspace method with shifts in the target set. This in turn yields approximations to eigenvalues of $F$. As the sizes of $(A,B)$ and, in turn, the cost of this approach increases with~$d$, there is clearly an advantage gained from the fact that, for a given tolerance~\texttt{reltol}, \alg often yields rational approximations of smaller, or otherwise at least comparable degree, compared to SV-AAA. This is in addition to the advantage of not requiring access to a split-form representation of~$F$.

 \subsection{Scattering problem}  \label{sec:scattering}
 
\begin{figure}[h]
\begin{minipage}[b]{0.49\linewidth}
\centering
\includegraphics[width=0.98\textwidth]{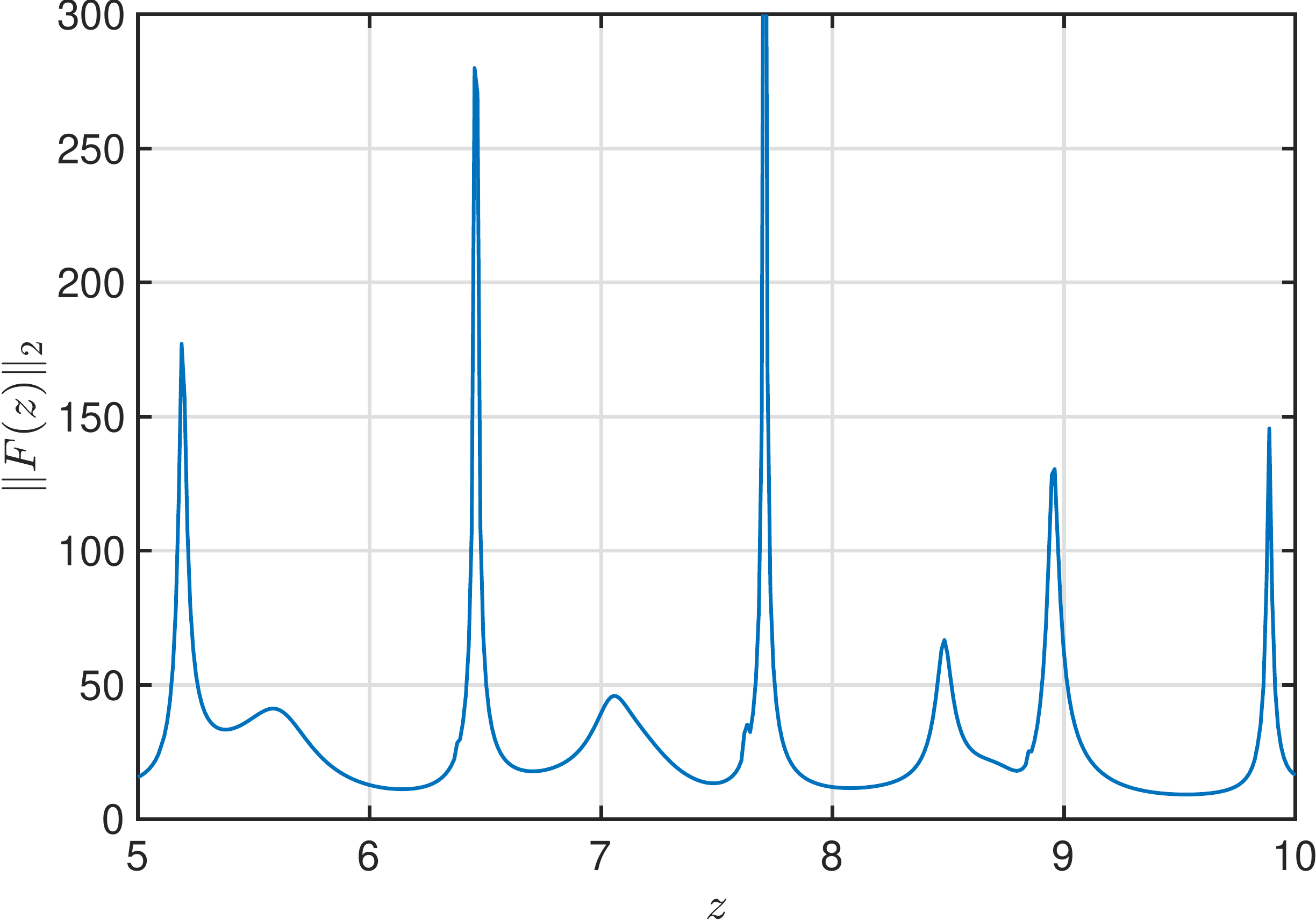}
\end{minipage}
\hfill
\begin{minipage}[b]{0.49\linewidth}
\centering
\includegraphics[width=\textwidth]{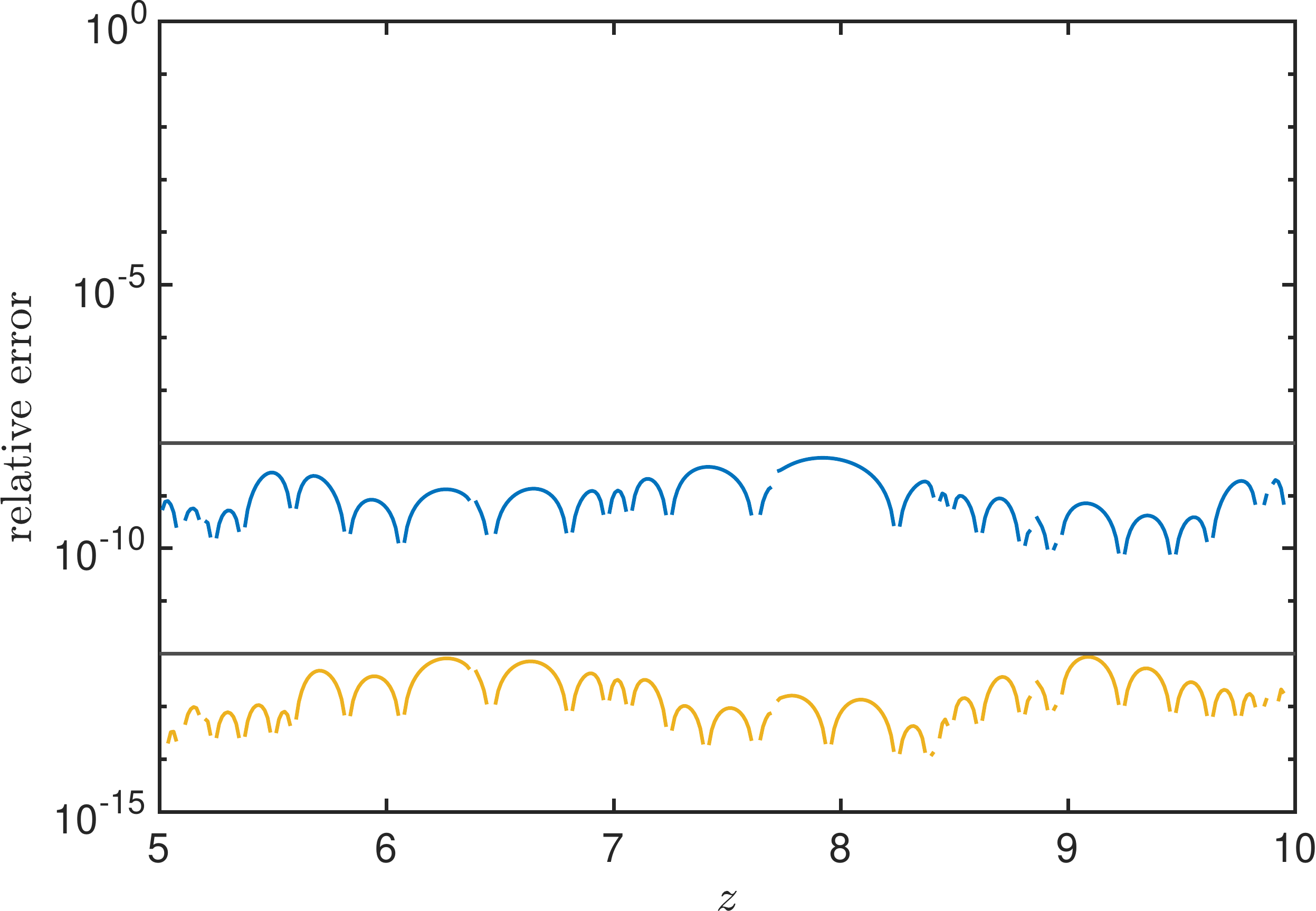}
\end{minipage}
\caption{Scattering problem. Left: the Euclidean norm of $\vf(z)$ as a function of $z$. Right: the relative error $\| \vf(z) - \vr^{(d)}(z) \|_{\infty} / \max_z \| \vf(z) \|_\infty$ for the rational approximants $\vr^{(d)}(z)$ with $\ell = 24$ probing vectors. The dark grey lines indicate the tolerances requested, namely, $10^{-8}$~and~$10^{-12}$.}
\label{fig:norm scattering}
\end{figure}

 We apply \alg to the Helmholtz equation with absorbing boundary conditions describing a scattering problem on the unit disc; see~\cite[Sec.~5.5.4]{Pradoverathesis}. The vector-valued function $\vf(z)$ containing the solution for a wavenumber $z>0$ is no longer given in split form and depends rationally on $z$:
 \[
  \vf(z) = (K - \mathrm{i} z C - z^2 M)^{-1} \mathbf{b}.
\]
Here, the stiffness matrix $K$, damping matrix $C$, and mass matrix $M$ are real non-symmetric sparse matrices of size 20054, obtained from a finite element discretization.  
The (complex) entries of $\vf(z)$ contain the nodal values of the finite element solution. 
The Euclidean norm of $\vf(z)$ for $z \in [5, 10]$ is depicted in Figure~\ref{fig:norm scattering}. Although there are no poles on the real axis, some are quite close to it, resulting in large peaks in $\|\vf(z)\|$. We therefore expect that a rather large degree for the rational approximant of $\vf$ will be needed.

The computational results of \alg applied to $\vf$ are reported in Table~\ref{tab:scattering}. The set $\Sigma$ contains 400 equidistant points in $[5,10]$. We observe that a large degree is indeed needed to get high accuracy. While the standard value of $\ell=4$ is performing decently, a larger sketch size is needed so that the error of the approximant is comparable to that of the surrogate. This behavior is reflected in our analysis in section~\ref{sec:fullanalysis}: according to Remark~\ref{remark:stablerank}, slower convergence of rational approximations signals larger stable rank, which in turn leads to less favorable probabilistic bounds in Theorem~\ref{thm:analysis} that are compensated by increasing $\ell$. However, let us stress that even for $\ell=24$ the rational approximation can be computed very quickly and it is still more than $500$ times faster than applying SV-AAA without sketching.
 
The timings in Table~\ref{tab:scattering} do not include the evaluation of $\vf$ and the error $\vr^{(d)}-\vf$ on the sampling set $\Sigma$, which is needed for all methods regardless of sketching. Since a large linear system has to be solved for each $z$, evaluating $\vf(z)$ is expensive. One of the benefits of rational approximation is that $\vr^{(d)}(z)$ can be evaluated much faster: the most accurate $\vr^{(d)}$ in Table~\ref{tab:scattering} can be evaluated in less than 0.002 seconds, whereas evaluating~$\vf$ requires 0.2 seconds.

 \begin{table}
\centering
\caption{Scattering problem from section~\ref{sec:scattering}: time to execute \alg given the evaluation of $\vf$ on~$\Sigma$, the  degree and relative error in the $\infty$-norm of the resulting approximation $\vr^{(d)}$ for 10~random realizations. The value $\ell=\infty$ indicates no sketching and corresponds to SV-AAA applied to the entire vector-valued function $\vf$.} \label{tab:scattering}
\begin{tabular}{lrccc} \toprule
\texttt{reltol} & $\ell$ & time (s.) &  degree & error \\ \midrule
1e-08  &  $\infty$  &  21.70 & 36.0  &  2.2e-09 \\  
      &  1  &  0.017 & 29.8  &  8.5e-06 \\  
      &  4  &  0.016 & 34.9  &  6.4e-08 \\  
      &  8  &  0.020 & 35.2  &  3.7e-08 \\  
      & 16  &  0.029 & 35.9  &  6.1e-09 \\  
      & 24  &  0.039 & 35.9  &  7.2e-09 \\ \midrule
1e-12 &  $\infty$  &  28.479 & 43.0  &  9.6e-13 \\  
      &  1  &  0.012 & 34.6  &  1.1e-07 \\  
      &  4  &  0.019 & 41.1  &  6.6e-11 \\  
      &  8  &  0.026 & 42.8  &  3.6e-12 \\  
      & 16  &  0.038 & 43.0  &  1.8e-12 \\  
      & 24  &  0.049 & 43.1  &  1.2e-12 \\     \bottomrule
 \end{tabular}
 \end{table}
 
\subsection{Boundary element matrices} \label{sec:bem}

Our last example is a nonlinear eigenvalue problem  that arises from the boundary integral formulation of an elliptic PDE eigenvalue problem~\cite{Steinbach2009}. More specifically, we consider the 2D Laplace eigenvalue problem on the Fichera corner  from~\cite[Section~4]{Effenberger2012a}. Applying the boundary element method (BEM) to this problem results in a matrix-valued function $F\colon \mathbb C \to \mathbb C^{n\times n}$ that is dense and not available in split form. Also, the entries of $F(z)$ are expensive to evaluate. Note, however, that hierarchical matrix techniques could be used to significantly accelerate assembly, resulting in a representation of $F(z)$ that allows for fast matrix-vector products~\cite{SauterSchwab}. Usually, the smallest (real) eigenvalues of $F$ are of interest. As the smallest eigenvalue is roughly $6.5$, we consider $z$ in the domain $[5,12]$, which is discretized by 200 equidistant points.

The number of boundary elements determines the size of the matrix $F(z) \in \C^{n \times n}$. We present two sets of numerical experiments depending on whether we can store the evaluations $F(z_i)$ for all sampling points $z_i \in \Sigma$ on our machine with 64 GB RAM.

\paragraph{Storage possible} The largest problem size that allows for storing all necessary evaluations of $F$ is $n=384$. The computational results are depicted in Table~\ref{tab:bem}. Like for the scattering problem, we see that larger sketch sizes are needed but \alg remains very fast and accurate. For example, with $\ell=16$  and for a stopping tolerance of 1e-08,  \alg is about 220 times faster than set-valued AAA applied to $\vecop(F(z))$ and achieves comparable accuracy. For a stopping tolerance of 1e-12,  \alg is about 350 times faster than set-valued AAA applied to $\vecop(F(z))$.

In all cases of  Table~\ref{tab:bem}, we excluded the 3.2~seconds it took to evaluate $F(z_i)$ for all $z_i \in \Sigma$. Both the set-valued AAA and \alg method require these evaluations.

 \begin{table}
\centering
 \caption{BEM problem from section~\ref{sec:bem}: time to execute \alg given the evaluation of $F$ on $\Sigma$, the  degree and relative error in the $\infty$-norm of the resulting \alg approximation $R^{(d)}$ for 10 random realizations. The value $\ell=\infty$ indicates no sketching and corresponds to SV-AAA applied to $\vecop(F(z))$. For $\ell>0$, full (non-tensorized) sketches where used.} \label{tab:bem}
\begin{tabular}{lcccc} \toprule
\texttt{reltol} & $\ell$ & time (s.) &  degree & error \\ \midrule
1e-08  &  $\infty$  &  15.147 & 11.0  &  1.4e-09 \\  
       &  1  &  0.023 & 7.4  &  1.9e-04 \\  
       &  4  &  0.030 & 10.0  &  1.3e-07 \\  
       &  8  &  0.045 & 10.3  &  2.6e-08 \\  
       & 16  &  0.067 & 10.8  &  5.2e-09 \\ \midrule
1e-12 &  $\infty$  &  21.737 & 14.0  &  2.0e-13 \\  
      &  1  &  0.023 & 9.1  &  2.1e-05 \\  
      &  4  &  0.034 & 12.9  &  5.7e-11 \\  
      &  8  &  0.046 & 13.3  &  7.1e-12 \\  
      & 16  &  0.061 & 14.0  &  2.6e-13 \\
\bottomrule
 \end{tabular}
 \end{table}
 
\paragraph{Storage not possible} For larger problems, $F(z_i)$ is evaluated when needed\footnote{A considerable part of the cost in assembling the BEM matrix $F(z_i)$ can be amortized when evaluating a few $z_i$ at the same time. We therefore evaluate and store $F$ in 20 values $z_1, \ldots, z_{20}$ at once and perform where possible all computations on $F(z_1), \ldots, F(z_{20})$ before moving on to the next set.} but never stored for all $z_i \in \Sigma$.  In Table~\ref{tab:bem_large} we list the results for tolerance $10^{-12}$. The degree and time for dense and tensor sketches are very similar and we therefore only show the results for tensor sketches. The errors for both the full and tensorized sketches are shown, even though they are also similar.

We observe that the runtime  of the whole algorithm is considerably higher. As expected it grows with~$N$. The degree and the error, on the other hand, remain very similar to those for the small problem and the main conclusion remains: for $\ell=1$ the \alg approximation is not accurate, but already for $\ell=4$ we obtain satisfactory results at the expense of only slightly increasing the runtime. In addition, even taking a large number of sketches $\ell=24$ is computationally feasible whereas applying SV-AAA to the original problem is far beyond what is possible on a normal desktop.

 \begin{table}
 \caption{BEM problem from section~\ref{sec:bem}: time to run \alg where $F(z_i)$ is evaluated on the fly, and the average degree and relative error in the $\infty$-norm of the resulting AAA approximation $R^{(d)}$ for 10 random realizations. 
 All numbers correspond to tensorized sketches except  the numbers in the column  \emph{error${}^*$} which reports the error obtained with full (non-tensorized) sketches. In all cases the tolerance \texttt{reltol} was $10^{-12}$.} \label{tab:bem_large}
\centering
\begin{tabular}{lccccc} \toprule
size $n$ &  $\ell$ & time (s.) &  degree & error & error${}^*$ \\ \midrule
384 & 1  &  26.134 & 9.4   &  1.2e-05 &  2.1e-05 \\  
    & 4  &  30.490 & 12.9  &  4.1e-11 &  5.7e-11 \\  
    & 8  &  31.132 & 13.4  &  8.7e-12 &  6.8e-12 \\  
    & 16  &  31.968 & 14.0  &  2.8e-13 &  2.6e-13 \\  
    & 24  &  32.017 & 14.0  &  2.2e-13 &  2.2e-13 \\ \midrule
864 & 1  &  130.172 & 9.1  &  7.6e-06   & 1.4e-05     \\  
    & 4  &  154.802 & 13.0  &  7.6e-11  & 5.0e-11     \\  
    & 8  &  159.993 & 13.8  &  6.6e-12  & 1.3e-11     \\  
    & 16  &  161.689 & 14.0  &  9.0e-13  & 7.5e-13     \\  
    & 24  &  161.614 & 14.0  &  6.9e-13  & 4.9e-13     \\ \midrule
1536  & 1  &  412.007 & 9.0  &  1.4e-05  & 8.6e-06 \\  
      & 4  &  489.553 & 12.9  &  5.5e-10 & 7.8e-11  \\  
      & 8  &  507.866 & 13.8  &  1.2e-11 & 1.5e-11  \\  
      & 16  &  512.706 & 14.0  &  1.4e-12 & 1.5e-12  \\  
      & 24  &  512.984 & 14.0  &  9.3e-13 & 1.2e-12 \\ \midrule
2400  & 1  &  1010.040 & 9.1   &  6.6e-06 & 1.5e-05 \\  
      & 4  &  1194.473 & 12.9  &  2.4e-10 & 1.3e-10 \\  
      & 8  &  1244.731 & 13.9  &  6.7e-12 & 1.2e-11 \\  
      & 16  &  1250.988 & 14.0  &  1.5e-12 & 1.8e-12 \\  
      & 24  &  1252.066 & 14.0  &  1.5e-12 & 1.4e-12 \\      
\bottomrule
 \end{tabular}
 \end{table}

\subsection{Empirical success probabilities}\label{sec: exp stats}

We verify numerically the bounds in the analysis from section~\ref{sec:analysis}. 
Let $\vf\colon \Omega \to \mathbb{K}^N$ be the vector-valued function to be approximated (e.g., the vectorization of $F(z)$ from a nonlinear eigenvalue problem). 
The aim is to investigate the reliability of the sketched estimator $$\mathrm{est}_\ell := \|  V_\ell^T \vf - \mathcal A( V_\ell^T \vf ) \|$$  in relation to the exact value $$\mathrm{ex} := \|  \vf - \mathcal A( \vf ) \|,$$
with $\| \cdot \|$ being the $L^2$ norm defined in~\eqref{eq:L2 norm}. Here, $\mathcal A$ is a rational approximant with fixed weights and fixed support points. More precisely,  we compute an AAA approximation of a certain degree for $\vf$ once, and then use its weights and support points to construct rational approximants for $V_\ell^T \vf$, one for each choice of $\ell$. The theoretical bounds in Theorem~\ref{thm:analysis} give us failure probabilities on the fidelity of $\mathrm{est}_\ell$, which imply the following estimates 
\begin{align*}
P_\mathrm{under} &= \mathrm{Prob}(  \mathrm{ex} / \tau < \mathrm{est}_\ell) = 1 - \mathrm{Prob}(\mathrm{ex}  \geq  \tau \, \mathrm{est}_\ell) \geq 1 - \mathrm{eq.~\eqref{eq:probunderestimate}} \\
P_\mathrm{over} &= \mathrm{Prob}(\mathrm{est}_\ell < \tau \,  \mathrm{ex}) = 1 - \mathrm{Prob}(\mathrm{est}_\ell \geq \tau \,  \mathrm{ex}) \geq 1 - \mathrm{eq.~\eqref{eq:proboverestimate}} \\
P_\mathrm{over \& under} &= \mathrm{Prob}(  \mathrm{ex} / \tau < \mathrm{est}_\ell < \tau \, \mathrm{ex}) \geq 1 - \mathrm{eq.~\eqref{eq:probunderestimate}} - \mathrm{eq.~\eqref{eq:proboverestimate}}
\end{align*}
for the success probabilities of under and/or over approximation of the exact value,~$\mathrm{ex}$, by a factor $\tau> 1$.

For evaluating the bounds above, we need the stable rank $\rho(\vh)$ of an $L^2$ function~$\vh$ defined on $\Omega$. We estimate this quantity by taking the corresponding stable rank $\rho(H)=\|H\|_F^2 / \|H\|_2^2$ of the matrix $H$ whose columns correspond to $\vh(z)$ with $z \in Z$ for a sufficiently dense finite sampling set $\Sigma \subset \Omega$.

The numerical results are shown in Table~\ref{tab:emp_failure}, with $\vf$ arising from the \texttt{buckling\_plate} problem that we already used in section~\ref{sec:nlvep experiments}. For each value of $\tau$ and $\ell$, the empirical success probabilities were calculated by taking $10^5$ random samples for $V_\ell$. The operator $\mathcal{A}$ was fixed throughout the experiment and corresponds to an AAA approximant of degree 25 with error $10^{-9}$. 

The table shows that the theoretical bounds are valid lower bounds when comparing to the empirical values. For larger values of $\ell$, there is a good correspondence. From a practical point of view, we see that when $\ell\geq 8$, the estimator is practically always within a tolerance of factor $\tau=3$. On the other hand, for $\ell=1$  the failure probability is not negligible: the estimator might fail $8\%$ of the time even with $\tau=10$.

\begin{table}
 \caption{Success probabilities of over and/or under approximation (empirical with $10^5$ samples and bounds according from Thm.~\ref{thm:analysis}) for the \texttt{buckling\_plate} problem when evaluating a fixed AAA rational approximant of degree 25. The stable rank satisfies $\rho \approx 1.05$.} \label{tab:emp_failure}
\centering
\begin{tabular}{p{0.02\textwidth}p{0.03\textwidth}p{0.1\textwidth}p{0.1\textwidth}p{0.1\textwidth}p{0.1\textwidth}p{0.1\textwidth}p{0.1\textwidth}} \toprule
  & &  \multicolumn{2}{c}{$P_\mathrm{over}$}  & \multicolumn{2}{c}{$P_\mathrm{under}$ } & \multicolumn{2}{c}{$P_\mathrm{over \& under}$ } \\
$\tau$  &  $\ell$ & empirical &  \emph{bound} & empirical & \emph{bound} & empirical &  \emph{bound} \\ \midrule
2 & 1 & 0.660 & \emph{0.607} & 0.959 & \emph{0.410} & 0.619 & \emph{0.018} \\  
  & 2 & 0.816 & \emph{0.768} & 0.985 & \emph{0.652} & 0.800 & \emph{0.420} \\  
  & 4 & 0.936 & \emph{0.901} & 0.998 & \emph{0.879} & 0.934 & \emph{0.780} \\  
  & 8 & 0.990 & \emph{0.977} & 1.000 & \emph{0.985} & 0.990 & \emph{0.963} \\  
  &16 & 1.000 & \emph{0.998} & 1.000 & \emph{1.000} & 1.000 & \emph{0.998} \\   \midrule
3 & 1 & 0.812 & \emph{0.732} & 0.998 & \emph{0.879} & 0.810 & \emph{0.611} \\  
  & 2 & 0.937 & \emph{0.889} & 1.000 & \emph{0.985} & 0.937 & \emph{0.875} \\  
  & 4 & 0.992 & \emph{0.976} & 1.000 & \emph{1.000} & 0.992 & \emph{0.976} \\  
  & 8 & 1.000 & \emph{0.999} & 1.000 & \emph{1.000} & 1.000 & \emph{0.999} \\   
  &16 & 1.000 & \emph{1.000} & 1.000 & \emph{1.000} & 1.000 & \emph{1.000} \\  \midrule
5 & 1 & 0.922 & \emph{0.837} & 1.000 & \emph{1.000} & 0.922 & \emph{0.837} \\  
  & 2 & 0.989 & \emph{0.959} & 1.000 & \emph{1.000} & 0.989 & \emph{0.959} \\  
  & 4 & 1.000 & \emph{0.997} & 1.000 & \emph{1.000} & 1.000 & \emph{0.997} \\  
  & 8 & 1.000 & \emph{1.000} & 1.000 & \emph{1.000} & 1.000 & \emph{1.000} \\  
  &16 & 1.000 & \emph{1.000} & 1.000 & \emph{1.000} & 1.000 & \emph{1.000} \\   \midrule
10& 1 & 0.978 & \emph{0.918} & 1.000 & \emph{1.000} & 0.978 & \emph{0.918} \\  
  & 2 & 0.999 & \emph{0.989} & 1.000 & \emph{1.000} & 0.999 & \emph{0.989} \\  
  & 4 & 1.000 & \emph{1.000} & 1.000 & \emph{1.000} & 1.000 & \emph{1.000} \\  
  & 8 & 1.000 & \emph{1.000} & 1.000 & \emph{1.000} & 1.000 & \emph{1.000} \\  
  &16 & 1.000 & \emph{1.000} & 1.000 & \emph{1.000} & 1.000 & \emph{1.000} \\ 
\bottomrule
 \end{tabular}
 \end{table}

\section{Conclusions} We have presented and analyzed a new randomized sketching approach which allows for the fast and efficient rational approximation of large-scale vector- and matrix-valued functions. Compared to the original surrogate-AAA approach in \cite{EG19}, our method \alg reliably achieves high approximation accuracy by using multiple (tensorized or non-tensorized) probing vectors. We have demonstrated the method's performance on a number of nonlinear functions arising in several applications. Compared to the set-valued AAA method from \cite{lietaert2018automatic}, our method works efficiently in the case when the split form of the function to approximate has a large number of terms, and even when the problem is only accessible via function evaluations. We believe that \alg is the first rational approximation method that combines these advantages. 

While our focus was on AAA and NLEVPs, let us highlight once more that our sketching approach is not limited to such settings. In principle, any linear (e.g., polynomial or fixed-denominator rational) approximation scheme applied to a large-scale vector-valued function can be accelerated  by this approach. For example, current efforts are underway to develop a sketched RKFIT method \cite{berljafa2017rkfit} as a replacement of the surrogate-AAA eigensolver in the MATLAB Rational Krylov Toolbox (the latter of which is currently using only $\ell=1$ sketching vector). 

We have used Gaussian sketches throughout this paper while subsampled randomized trigonometric transforms and sparse sketching operators are also popular; see~\cite{murray2023randomized} and the references therein. We decided not to use those in our experiments as (i) even with full Gaussian sketches, the sketching runtime is relatively small compared to that for the computation of the rational approximants, and (ii) the available theory is much less developed for non-Gaussian sketching operators. In particular, the \emph{small-sample} estimates of Theorem~\ref{thm:analysis} are very different from the oblivious subspace embedding (OSE) property usually established for non-Gaussian sketches. Our analysis implies low failure property as~$\tau$ increases. Such a property (sometimes referred to as small ball probability) is not implied by OSE and has, to the best of our knowledge, not been proven for non-Gaussian sketches---even for the case of sketching matrices instead of functions.

There are a number of interesting research directions arising from this work. This includes a possible extension to the multivariate case. A multivariate p-AAA algorithm has recently been proposed in \cite{rodriguez2020p} but it is not immediately clear whether the sketching idea pursued here can be extended to this algorithm. Another potential improvement of \alg in the case of many sampling points is to replace the SVDs for the least-squares problems \eqref{eq:min L w} by another sketching-based least squares solver such as~\cite{rokhlin2008fast}, similarly to what has been done in~\cite{nakatsukasa2022fast}. Finally, we hope that the analysis provided in this paper might shed some more light onto the accuracy of contour integral-based solvers for linear eigenvalue problems $A\vx=\lambda \vx$. These methods can be viewed as pole finders for the resolvent $(A-zI)^{-1}$ after random tensorized probing of the form $\mathbf{u}^T (A-zI)^{-1} \mathbf{v}$.

\section*{Acknowledgments}
We thank Davide Pradovera for providing us with the code for the scattering problem considered in section~\ref{sec:scattering}. We also thank the two anonymous referees who have provided useful comments and critical insights. B.V.~was supported by the Swiss National Science Foundation (grant number 192129).

\bibliographystyle{siam} 
\bibliography{refs}

\end{document}